\documentclass{article}
\usepackage{amsmath,amscd,amsthm}
\begin{document}
\author{Oliver Attie and Sylvain Cappell}
\title{Bott Integrability and Higher Bott Integrability; Higher Cheeger-Simons and Godbillon-Vey Invariants}
\maketitle
\begin{abstract}
This paper studies the interaction of $\pi_1(M)$ for a $C^\infty$ manifold
$M$  with Bott's original 
obstruction to integrability, and with differential geometric invariants
such as Godbillon-Vey and Cheeger-Simons invariants of a foliation.
We prove that the ring of higher Pontrjagin and higher Chern classes of 
an integrable subbundle $E$ of the tangent bundle of a manifold vanishes above
dimension $2k$ where $k=dim(TM/E)$, and where the higher Pontrjagin and 
Chern rings are rings generated by $i^*y \cup p_j(TM/E)$ and by
$i^*y \cup c_j(TM/E)$ respectively, with $p_j$ the $j$-th Pontrjagin class,
$c_j$ the $j$-th Chern class, $i:M \to B\pi$ and $\pi=\pi_1(BG)$, 
where $BG$ is the classifying space of the holonomy groupoid corresponding to 
$E$ and $y \in H^*(B\pi)$, provided that the fundamental group of $BG$ 
satisfies the Novikov conjecture. In addition, we show the vanishing of
higher Pontrjagin and Chern rings generated by $i^*x \cup p_j(TM/E)$, and
by $i^*x \cup c_j(TM/E)$ as before but with $i:M \to BG$, $BG$ as above
and $x \in H^*(BG)$ provided $(M,\mathcal{F})$ satisfied the foliated
Novikov conjecture, where $\mathcal{F}$ is the foliation whose tangent bundle
is $E$. We give examples of this obstruction
and of higher Godbillon-Vey and Cheeger-Simons invariants. 
\end{abstract}
\section{Introduction} If $M$ is a $C^\infty$ manifold, a subbundle
$E \subset TM$ is said to be integrable if its space of sections is 
closed under Lie bracket, i.e., $X,Y \in \Gamma(E)$, implies $[X,Y] \in 
\Gamma(E)$. 
\newtheorem{thm}{Theorem}[section]
\theoremstyle{corollary}
\newtheorem*{corollary}{Corollary 4.2}
\newtheorem*{corollary1}{Corollary 4.4}
\begin{thm}[Frobenius] \cite{CC} If $E$ is a subbundle of the tangent bundle, 
then $E$ is integrable if and only if $E=T(\mathcal{F})$ for a foliation
$\mathcal{F}$ of $M$.
\end{thm}

In 1966 S.S.Chern \cite{Chern} asked whether every subbundle
of the tangent bundle of a compact manifold is deformable to an integrable
one. Raoul Bott answered
this question in 1968, showing that there indeed was an obstruction 
involving the vanishing in a certain range of dimensions of  characterstic
cohomology classes \cite{Bott}. Bott's integrability 
condition was proven by him in two ways. First, by using differential 
geometry \cite{Bott}, constructing a connection (called a basic connection)
which is flat along the 
leaves of the foliation. Herbert Shulman, in his 1972 U.C.Berkeley Ph.D. thesis,
rephrased Bott integrability as a statement about the terms in the
spectral sequence associated to a semi-simplicial model for the foliated
manifold.
Then \cite{BSS} in 1976, Bott used
algebraic topology, the cohomology of the Haefliger space, to prove his 
integrability condition. A similar proof was obtained for foliated bundles
by Kamber and Tondeur in 1975 \cite{KT}. One consequence of the 
integrability theorem is the existence of characteristic classes for 
foliations \cite{BC,CS, KT}.

In this paper, we begin by presenting a new formulation of the proof of Bott's
Integrability Theorem:
\newtheorem*{theorem}{Theorem 3.2}
\begin{theorem}[Bott] Let $Pont^*(TM/E)$ be the Pontrjagin ring of $TM/E$.
  Then if $E$ is integrable, for $q>2k$, $Pont^q(TM/E)=0$.
  \end{theorem}
Our proof uses, as did Bott's, a result by Shulman on characteristic
differential forms but in place of the basic connection theory we use an
assembly map in surgery theory. This approach is then used to generalize
the theorem to obtain higher obstructions to integrability in some cases of
nonsimply connected manifolds and to give examples of these new obstructions,
a first improvement on Bott's integrability theorem since 1968.
The higher characteristic classes that appear in the new obstructions,
are certain products of characteristic classes with natural cohomology
classes, e.g., those coming from the fundamental group or 
the classifying space of the holonomy groupoid. Such classes play a role in the
classical Novikov conjecture about such characteristic classes being
oriented homotopy invariants. This conjecture has been verified for many
large and important classes of fundamental groups.
\newtheorem*{theorem1}{Theorem 4.1}
\begin{theorem1} Let $M$ be a compact manifold of dimension $n$, suppose
  $E \subset TM$ is a subbundle of the tangent bundle of $M$. Let $E$ be
  integrable, $BG$ the classifying space of the foliation groupoid of the
  corresponding foliation and $\pi=\pi_1(BG)$. Suppose the Novikov conjecture
  is true for $\pi$. Then the higher Pontrjagin ring $Pont^q(\pi)(TM/E)$
  vanishes above dimension $2k$, where $k$ is the dimension of $TM/E$.
\end{theorem1}
This theorem also  holds for classes
in the cohomology of the classifying space of the holonomy groupoid
 provided that the Novikov Conjecture for foliations holds
for that holonomy groupoid \cite{AttCap}, as we shall see in section 6.

We can apply to the hypothesis of Theorem 4.1 the following result for
$(k-1)=1$.
\begin{thm}[Haefliger \cite{Haefliger}] Let $\mathcal{F}$ be a foliation
  on a manifold $X$ so that the holonomy coverings of the leaves are
  $(k-1)$-connected. Then the holonomy groupoid $G$ of $\mathcal{F}$
  considered as a $G$-principal bundle with base $X$ by the end projection
  is $k$-universal. Thus the space $X$ itself is $k$-classifying and the
  map $i$ of $X$ to $BG$ is $k$-connected.
\end{thm}
That then yields:
\begin{corollary} Suppose the leaves of the foliation $\mathcal{F}$ on $M$
  above are all simply connected. Then $\pi_1(M)=\pi_1(BG)=\pi$. Suppose
  the Novikov conjecture is true for $\pi$. Then the higher Pontrjagin ring
  of $TM/E$, $Pont^q(\pi)(TM/E)$ vanishes for $q>2k$ where $k=dim(TM/E)$.
  \end{corollary}
We also have results that do not depend on the Novikov conjecture. One such
result is:
\begin{corollary1} Let $M$ be a closed, oriented, smooth manifold with an
  integrable subbundle $E$ of the tangent bundle $TM$ of $M$,
  and a homomorphism $f:\pi_1(M)\to \textbf{Z}$ so that the
  corresponding rational cohomology class $f^*(j)$, $j$ the generator
  of $H^1(\textbf{Z})$, pulls back to 0 on each leaf of the foliation.
  More generally, suppose that $\gamma \in H^*(M;\textbf{Q})$ belongs to the
  cohomology algebra generated by $H^1(M;\textbf{Q})$ and pulls back
  to 0 on each leaf of the foliation. Then the higher signature class
  $f^*(j)\cup L_q(TM/E)$ or $\gamma \cup L_q(TM/E)$ vanishes for $q+1$
  resp. $q+dim(\gamma)$ larger than $2dim(TM/E)$.
  \end{corollary1}

The leaves of a foliation are smooth manifolds of bounded geometry. We 
can therefore use the algebraic topology \cite{Att,Engel,J,Roe,Spakula}
of manifolds 
of bounded geometry to study foliations, and derive the Bott integrability 
theorem. We use the assembly map in surgery theory for manifolds of bounded
geometry to show that since Bott's integrability condition holds for the
normal bundles of the leaves, considered as vector bundles of bounded geometry
and having a Pontrjagin ring for the $L^\infty$ Pontrjagin classes of
the normal bundle, it holds for the normal bundle to the foliation, to which
they assemble. (In a related paper, Baum and Connes \cite{BaumConnes} state a
Novikov Conjecture for foliations using the K-theory of operator algebras
and Connes' foliation algebra).
We carry out the assembly by using a cosheaf of L-homology
groups to detect the Pontrjagin classes.

We use this integrability theorem to construct new characteristic classes
in Section 7:
a higher Godbillon-Vey class, a higher Cheeger-Simons class,
and a higher Maslov index in Example 7.4.
Classes beyond Godbillon-Vey and its higher codimension versions have been
defined by S.Hurder and D.McDuff \cite{Hurder3, McDuff}. Connes
\cite{Connes2} has proven the leafwise homotopy invariance of the L-class
evaluated on secondary classes, and that secondary classes are obstructions
to positive scalar curvature of the tangent bundle to the foliation.

Bott's obstruction is analogous to the obstruction to fibering over a manifold
\cite{Quinn}, where the leaves of the foliation replace the compact 
manifold fiber. As a result of this approach we get a new set of obstructions, 
the ring of higher Pontrjagin classes, which must also vanish for a subbundle 
to be integrable. These classes yield higher obstructions
 only when the manifold $M$ has infinite 
fundamental group. Many of the constructions go through easily in the smooth 
case. We will further discuss the analogy with the blocked surgery of
\cite{Quinn} in \cite{AttCap}, where we will construct an analogous blocked
surgery exact sequence for foliations, which is independent of this paper.

The authors wish to acknowledge our dear late friend Andrew Ranicki for the 
definition of uniformly finite L-homology. We also wish to acknowledge
Steve Hurder and Shmuel Weinberger for helpful comments on this paper.
\section{Definitions}
The present melding of foliation theory with surgery-theoretic ideas
applied to manifolds of bounded geometry arising as leaves of foliations
and also some aspects of operator theory which we develop here
(and further in \cite{AttCap}) involves a list
of definitions or adaptations as follows.
\newtheorem{dfn}{Definition}[section]
\begin{dfn} A simplicial complex has bounded geometry if there is a uniform
bound on the number of simplices in the link of each vertex \cite{Att}.
\end{dfn}
\begin{dfn} A complete paracompact Riemannian manifold $M$ is said to have 
bounded geometry if its injectivity radius $inj_M > c > 0$ for some constant
$c$ and its sectional curvature is bounded in absolute value.
\end{dfn}
\newtheorem{rk}{Remark}[section]
\begin{rk}The leaves of a foliation of a compact manifold are manifolds 
of bounded geometry, and conversely \cite{AB} every manifold of bounded 
geometry is the leaf of a compact foliated space.
But it is not true that every manifold of bounded geometry is the leaf
of a $C^1$ foliation of a compact manifold \cite{AH}.
\end{rk}
\begin{thm} Every smooth manifold $M$ of bounded geometry has a triangulation
as a simplicial complex of bounded geometry. Conversely, if $M$ is a 
simplicial complex of bounded geometry which is a triangulation of a smooth
manifold, then this smooth manifold admits a metric of bounded geometry
with respect to which it is quasi-isometric to $M$ \cite{Att}.
\end{thm}
\begin{dfn} Let $X$ be a simplicial complex of bounded geometry. 
The uniformly finite homology groups of $X$ with coefficients in an abelian normed 
group $(G, \mid \cdot \mid)$ denoted $H^{uff}_*(X; G, \mid \cdot \mid)$
are defined to be the homology groups of the complex of infinite simplicial
chains whose coefficients are in $l^\infty$ with respect to the norm
$\mid \cdot \mid$ on $G$ \cite{Att}.
\end{dfn}
\begin{dfn} Denote by $\Omega^p_\beta(M)$ the Banach space of $p$-forms on a
complete, oriented Riemannian manifold $M$ which are bounded in the norm
$$\parallel \alpha \parallel=sup\{\mid \alpha(x)\mid+\mid d\alpha(x)\mid
: x \in M\}.$$
This gives rise to a complex $d_i:\Omega_\beta^i(M)\to \Omega^{i+1}_\beta(M).$
The bounded de Rham groups are defined by
$$H_\beta^p(M)=[\mbox{Ker }d_p]/[\mbox{Im }d_{p-1}],$$
where we are not taking the closure of $\mbox{Im }d$ in this 
definition \cite{Att}.
\end{dfn}
\begin{dfn} The $L^\infty$ Pontrjagin classes are defined (by Januszkiewicz)
 to be characteristic classes in $H^{4*}_\beta(M)$, the $L^\infty$ de Rham 
cohomology of $M$, where $M$ is a manifold of bounded geometry. These are 
the result of a Chern-Weil homomorphism defined by Januszkiewicz \cite{J}.
\end{dfn}
\begin{dfn} The $L^\infty$ Chern classes are defined (by Januszkiewicz) to be
  characteristic classes in $H^{2*}_\beta(M)$, the $L^\infty$ de Rham
  cohomology of $M$, where $M$ is a manifold of bounded geometry. These are
  the result of the Chern-Weil homomorphism defined by Januszkiewicz \cite{J}.
\end{dfn}
\begin{dfn} Let $M$ be a manifold of bounded geometry. 
If $E \to M$ is a vector bundle, $\nabla$ a connection on $E$,
then $(E,\nabla)$ has bounded geometry if the curvature tensor of $E$  
and all its covariant derivatives are bounded \cite{Roe}. 
The direct sum of such $(E, \nabla)$ and 
$(E^\prime, \nabla^\prime)$ is the vector bundle of bounded geometry 
$(E \oplus E^\prime, \nabla \oplus \nabla^\prime)$.
\end{dfn}
\begin{dfn} The $L^\infty$ Hirzebruch L-classes are defined in terms of the 
$L^\infty$ Pontrjagin classes as the polynomials associated to the genus 
$\mathcal{L}$ given by $S(\mathcal{L})=\sqrt{x}/tanh\sqrt{x}$. They are 
  denoted $L_0^{(\infty)}(\xi),L_1^{(\infty)}(\xi),...,L_k^{(\infty)}(\xi)$,
  for a given vector bundle $\xi$ of bounded geometry.
\end{dfn}
\begin{thm}[Roe] Let $M$ be a manifold of bounded geometry, $E \to M$ a 
vector bundle. Then the following are equivalent: 
\par\noindent
i.$E \to M$ is a vector bundle whose curvature tensor and all its derivatives
are uniformly bounded.
\par\noindent
ii. The Christoffel symbols of $E$ are bounded, as are all their derivatives, 
these bounds being independent of $x\in M$.
\par\noindent 
iii. The transition functions of $E$ along with all their derivatives are 
bounded, bounds being the same for all transition functions \cite{Roe, Shubin}.
\end{thm} 
\begin{dfn} Let $X$ be a manifold of bounded geometry. 
$K_{uf}^*(X)$ is the Grothendieck group of complex vector bundles of bounded 
geometry with their connections under direct sum. $KO_{uf}^*(X)$ is the 
Grothendieck group of real vector bundles of bounded geometry with their 
connections under direct sum \cite{Engel}, which we call the uniformly
finite KO-theory of $X$.
\end{dfn}
\begin{dfn} \cite{Engel} A graded Hilbert space is a Hilbert space $H$ with
  a decomposition $H=H^+ \oplus H^-$ into closed orthogonal spaces. This
  is the same as prescribing a grading operator $\epsilon$ whose
  $\pm 1$-eigenspaces are $H^\pm$ and such that $\epsilon$ is selfadjoint
  and unitary. If $H$ is a graded space, then its opposite is the graded
  space $H^{op}$ whose underlying vector space is $H$, but with reversed
  grading, i.e., $(H^{op})^+=H^-$ and $(H^{op})^-=H^+$. An operator on
  a graded space $H$ is called even if it maps $H^\pm$ to $H^\pm$ and it
  is called odd if it maps $H^\pm$ to $H^\mp$. Equivalently an operator is
  even if it commutes with the grading operator $\epsilon$ of $H$, and it
  is odd if it anti-commutes with it.
  \end{dfn}
\begin{dfn} \cite{Engel} Let $p$ be a non-negative integer.
  A $p$-multigraded Hilbert
  space is a graded Hilbert space which is equipped with $p$
  operators, which are odd and unitary $\epsilon_1,...,\epsilon_p$ so that
  $\epsilon_i\epsilon_j+\epsilon_j\epsilon_i=0$ for $i \ne j$ and
  $\epsilon_j^2=-1$ for all $j$.
  If $H$ is a $p$-multigraded Hilbert space, then an operator $H$ is called
  multigraded if it commutes with the multigraded operators $\epsilon_1, ...,
  \epsilon_p$ of $H$.
\end{dfn}
\begin{dfn}\cite{Engel} Let p be an integer greater than -1. A triple $(H,\rho,T)$
  consisting of
  \begin{enumerate}
  \item [i.] A separable $p$-multigraded Hilbert space $H$.
  \item [ii.] A representation $\rho:C_0(X) \to \mathcal{B}(H)$ by even,
    multigraded operators
  \item[iii.] An odd multigraded operator $T \in \mathcal{B}(H)$ so that
    the operators $T^2-1$ and $T-T^*$ are locally compact and the operator
    $T$ itself is pseudocal.
  \end{enumerate}
  is called a $p$-multigraded Fredholm module over $X$.
\end{dfn}
\begin{dfn}\cite{Engel} A collection of operators $\mathcal{A}\subset \mathcal{K}(L^2(E))$
  is said to be uniformly approximable, if for every $\epsilon > 0$ there
  is an $N>0$ so that for every $T \in \mathcal{A}$ there is a rank-$N$
  operator $k$ with $\parallel T-k \parallel < \epsilon$.
\end{dfn}
\begin{dfn}\cite{Engel}
  $$L-Lip_R(X)=\{f\in C_c(X)\mid f\mbox{ is L-Lipschitz }, diam(supp(f))\le R
  \mbox{ and }\parallel f \parallel_\infty \le 1\}$$
\end{dfn}
\begin{dfn}\cite{Engel} Let $T \in \mathcal{B}(H)$ be an operator on a Hilbert space $H$
  and $\rho:C_0(X)\to \mathcal{B}(H)$ a representation. We say that $T$ is
  unifomly locally compact, if for every $R,L>0$ the collection
  $$\{\rho(f)T,T\rho(f) \mid f \in L-Lip_R(X)\}$$
  is uniformly approximable. We say that $T$ is uniformly pseudolocal, if
  for every $R,L>0$ the collection
  $$\{[T,\rho(f)] \mid f \in L-Lip_R(X)\}$$
  is uniformly approximable.
  \end{dfn}
\begin{dfn}\cite{Engel} A Fredholm module $(H,\rho,T)$ is uniform, if $T$ is uniformly
  pseudolocal and the operators $T^2-1$ and $T-T^*$ are uniformly locally
  compact.
  \end{dfn}
\begin{dfn}\cite{Engel} We define the analytic uniformly finite K-homology \break
  $K_p^{anal,uf}(X)$
  of a locally compact and separable metric space $X$ to be the abelian group
  generated by unitary equivalence classes of $p$-multigraded uniform Fredholm
  modules with the relations
  \begin{enumerate}
  \item[i.] If $x$ and $y$ are operator homotopic, then $[x]=[y]$.
  \item[ii.] $[x]+[y]=[x \oplus y]$
  \end{enumerate}
  where $x$ and $y$ are $p$-multigraded uniform Fredholm modules.
\end{dfn}
\begin{thm}[Engel] The uniformly finite K-theory of $X$ is isomorphic to
  the operator K-theory of $C_u(X)$, the $C^*$-algebra of all bounded,
  uniformly continuous functions on $X$:
  $$K_{uf}^p(X)=K_p(C_u(X))$$
  \end{thm}
\begin{dfn}\cite{Engel} Engel defines a cap product:
  $$\cap:K_{uf}^p(X)\otimes K_q^{anal,uf}(X)\to K_{q-p}^{anal,uf}(X)$$
  as follows:
\par\noindent
Let $P$ be a projection in $Mat_{n \times n}(C_u(X))$, where $C_u(X)$ is
the $C^*$-algebra of all bounded, uniformly continuous functions on $X$,
and let
  $(H,\rho,T)$ be a uniform Fredholm module. We set $H_n=H\otimes \textbf{C}^n$,
  $\rho_n(-)=\rho(-)\otimes id_{\textbf{C}^n}$, $T_n=T\otimes id_{\textbf{C}^n}$
  and by $\pi$ we denote the matrix
  $\pi_{ij}=\rho(P_{ij}) \in Mat_{n\times n}(\mathcal{B}(H))=\mathcal{B}(H_n)$.
  Then $(\pi H_n, \pi \rho_n \pi, \pi T_n \pi)$ is a uniform Fredholm
  module, with an induced (multi-)grading if $(H,\rho,T)$ was (multi-)graded.
  This construction is compatible with relations defining K-theory and
  uniform K-homology and is bilinear.
\end{dfn}

\begin{thm}[Engel]\cite{Engel} Let $M$ be an $m$-dimensional $KU$-oriented manifold of
  bounded geometry without boundary. Then the cap product
  $$\cap [M]:K_{uf}^*(M) \to K_{m-*}^{anal,uf}(M)$$
  with its uniform fundamental class $[M]\in K_m^{anal,uf}(M)$ is an isomorphism.
\end{thm}
\begin{thm}\cite{Engel} Let $M$ be an $m$-dimensional $KO$-oriented manifold of bounded geometry
  without boundary. Then the cap product
  $$\cap [M]:KO_{uf}^*(M) \to KO_{m-*}^{anal,uf}(M)$$
  with its uniform fundamental class $[M]\in KO_m^{anal,uf}(M)$ is an
  isomorphism.
  \end{thm}
  \begin{dfn}\cite{BHSchick} Let $X$ be simplicial complex of bounded geometry, $Y$ a closed
  subcomplex of $X$. A $K$-cycle for the pair $(X,Y)$ is a triple $(M,E,\phi)$
  consisting of
  \begin{enumerate}
  \item [i.] A smooth manifold $M$ of bounded geometry equipped with a $Spin^c$ structure.
  \item [ii.]A smooth Hermitian $bg$ vector bundle $(E,\nabla)$ on $M$.
  \item [iii.] A continuous map $\phi:M \to X$ such that $\phi[\partial M]\subseteq Y$.
  \end{enumerate}
  Two uniformly finite  K-cycles are isomorphic if there are compatible
  isomorphisms of all of the
  above three components in the definition.
\end{dfn}
\begin{dfn}\cite{BHSchick} If $(M,E,\phi)$ and $(M^\prime, E^\prime, \phi^\prime)$ are two
  uniformly finite K-cycles for $(X,Y)$ then their disjoint union is
  the K-cycle $(M \cup M^\prime, E \cup E^\prime, \phi \cup \phi^\prime)$.
\end{dfn}
\begin{dfn}\cite{BHSchick} If $(M,E,\phi)$ is a K-cycle for $(X,Y)$, then its opposite is the
  K-cycle $(-M,E,\phi)$, where $-M$ denotes the manifold $M$ equipped with the
  opposite $Spin^c$-structure.
\end{dfn}
\begin{dfn}\cite{BHSchick} A bordism of K-cycles for the pair $(X,Y)$ consists of the
  following data:
  \begin{enumerate}
  \item[i.] A smooth manifold of bounded geometry $L$ equipped with a
    $Spin^c$-structure.
  \item[ii.] A smooth, Hermitian $bg$ vector bundle $(F,\nabla)$ over $L$.
  \item[iii.] A continuous map $\Phi:L \to X$.
  \item[iv.] A smooth map $f:\partial L \to \textbf{R}$ for which
    $\pm 1$ are regular values, and for which $\Phi[f^{-1}[-1,1]]\subseteq Y$.
  \end{enumerate}
\end{dfn}
\begin{dfn}\cite{BHSchick} Let $(M,E,\phi)$ be a K-cycle for $(X,Y)$ and let $W$ be a
  $Spin^c$-vector vector bundle over $M$ with even dimensional fibers. Let
  $Z$ be the sphere bundle of $W \oplus \textbf{1}$, where $\textbf{1}$ is
  trivial rank-one real vector bundle. The vertical tangent bundle of $Z$ has
  a natural $Spin^c$-structure. Denote by $S_V$ the corresponding reduced
  spinor bundle and let $F=S^*_{V,+}$. In other words, define $F$ to be
  the dual of the even-graded part of the $\textbf{Z}/2$-graded bundle
  $S_V$. The modification $(M,E,\phi)$ associated to $W$ is the K-cycle
  $(Z,F\otimes\pi^*E, \phi\circ\pi)$
\end{dfn}
\begin{dfn}\cite{BHSchick} Denote by $K^{uf}(X,Y)$ the set of equivalence classes of uniformly
  finite K-cycles over $(X,Y)$ for the equivalence relations:
  \begin{enumerate}
  \item[i.] If $(M,E_1,\phi)$ and $(M,E_2,\phi)$ are two uniformly finite
    K-cycles with the same $Spin^c$-manifold $M$ and map $\phi:M\to X$, then
    $$(M\cup M,E_1 \cup E_2, \phi \cup \phi)\sim (M,E_1 \oplus E_2, \phi)$$
  \item[ii.] If $(M_1, E_1, \phi_1)$ and $(M_2, E_2, \phi_2)$ are bordant
    K-cycles then
    $$(M_1, E_1, \phi_1) \sim (M_2, E_2, \phi_2)$$
  \item[iii.] If $(M,E,\phi)$ is a uniformly finite K-cycle, and if $W$ is an
    even-dimensional $Spin^c$-vector bundle over $M$, then
    $$(M,E,\phi) \sim (Z,F\otimes \pi^*E, \phi\circ\pi)$$
    where $(Z,F\otimes \pi^*E,\phi \circ \pi)$ is the previous modification of
    $(M,E,\phi)$.
  \end{enumerate}
  The set $K^{uf}(X,Y)$ is an abelian group. The addition operation is given
  by disjoint union and the additive inverse of a cycle is given by reversing
  the $Spin^c$-structure.
\end{dfn}
We can by considering real Spin-structures and real bounded geometry vector
bundles in the above definitions obtain $KO^{uf}(X,Y)$. The subgroups of
$K^{uf}(X,Y)$ consisting of cycles $(M,E,\phi)$ for which $M$ is even
dimensional is $K^{uf}_0(X,Y)$, the subgroups of $K^{uf}(X,Y)$ consisting
of cycles $(M,E,\phi)$ for which $M$ is odd dimensional is $K_1^{uf}(X,Y)$.
Similarly for $KO^{uf}(X,Y)$ the subgroup of cycles $(M,E,\phi)$ and $M$
is of dimension $i$ mod 8 is $KO_i^{uf}(X,Y)$. For more on this see
\cite{RSV}.
\begin{dfn} \cite{BHSchick} We associate to each uniformly finite K-cycle $(M,E,\phi)$ for
  $(X,Y)$ a class $<M,E,\phi>$ in Kasparov K-homology. Denote by $M^\circ$
  the interior of $M$, which is a $bg$ $Spin^c$-manifold. The
  $Spin^c$-structure on $M$
  determines a spinor bundle $S$ on $M^\circ$ by restriction, and of course
  the complex vector bundle $E$ also restricts to $M^\circ$. The tensor
  product $S\otimes E$ is a Dirac bundle over $M^\circ$ and if $D_E$ is
  an associated Dirac operator, we can form the class
  $$[D_E] \in K^{anal,uf}_n(M^\circ).$$
  The map $\phi:M \to X$ restricts to a proper map from $M^\circ$ to $X\backslash Y$ and
  we can therefore form the class
  $$\phi_*[D_E]\in K^{anal,uf}_n(X,Y)$$
\end{dfn}
\begin{thm} (The analogue for ordinary K-homology is proven in \cite{BHSchick}).
The correspondance $(M,E,\phi) \mapsto \phi_*[D_E]$ determines a
functorial isomorphism
$$\mu:K_*^{uf}(X,Y) \to K_*^{anal,uf}(X,Y)$$
\end{thm}
\begin{dfn} A Lipschitz homotopy equivalence $f:X \to Y$ with a
  Lipschitz map $g:Y \to X$ so that $f\circ g$ is Lipschitz homotopic
  to the identity and $g \circ f$ is Lipschitz homotopic to the identity.
  \end{dfn}
\begin{dfn} A spectrum with a metric space structure $E$ is a spectrum
  $E=\{E_n\}$, with a metric space structure on $E_n$ and Lipschitz
  homotopy equivalences $f_n:E_n \to \Omega E_{n+1}$ with Lipschitz
  constant of the composition $f_n\circ ... \circ f_{n+k}:E_n \to
  \Omega^k E_{n+k}$ is uniformly bounded for all $k$, and so that
  $\underset{k\to\infty}{\lim}f_n\circ...\circ f_{n+k}$ is a Lipschitz homotopy
  equivalence. 
\end{dfn}
\begin{rk} This allows us to define $QX$ for $X=E_0$ the 0-space of the
  spectrum $E$ where $QX=\Omega^\infty\Sigma^\infty X$.
  \end{rk}
\begin{dfn} The uniformly finite cohomology $E_{uf}^*(X)$, where $E$ is a
  CW-spectrum equipped with a fixed structure as a metric space and $X$ is a
  simplicial complex is defined as
  $$E_{uf}^*(X)=[X:\Omega^*(E)]^{Lip}$$
  where $[X:Y]^{Lip}$ is Lipschitz homotopy classes of Lipschitz maps with
  uniform Lipschitz constant and $\Omega^*(E)$ is given the metric space
  structure induced from $E$, by considering $\Omega(X)$ to be a function
  space with the $L^\infty$ metric. This is not independent of the topological
  type of $E$ or its metric. For example let $E_0=\underset{\to}{\lim}(Symm^k(S^n))$ where
  $Symm^k$ is the $k$-th symmetric power. Then the Dold-Thom theorem gives
  $E_0=K(\textbf{Z}, n)$. On the other hand $X$ could be a CW complex with
  finite skeleta the skeleta of $K(\textbf{Z},n)$. The latter will be
  assumed to represent uniformly finite cohomology in this paper.
  \end{dfn}
\begin{dfn} Let $M$ be a manifold of bounded geometry, $\alpha$ a vector
  bundle of bounded geometry over $M$. The Thom space $M^\alpha$ is defined
  to be $D(\alpha)/S(\alpha)$ where $D(\alpha),S(\alpha)$ are the unit disk
  and sphere bundles of $\alpha$, considered as a simplicial complex of
  bounded geometry. Then if $M$ is a manifold of bounded geometry
  quasi-isometrically embedded into Euclidean space $\textbf{R}^s$ with
  normal bundle $\nu$, we define the S-dual of $M \cup pt$ as $M^\nu$.
\end{dfn}
\begin{dfn} \cite{AM,Att} Let $X$ be a space controlled over a metric space $Z$ by a
  control map $p$. Denote by $\mathcal{P}$ the category of metric balls
  in $Z$ with morphisms given by inclusions. Define $\mathcal{P}G_1(X)$
  to be the category whose objects are pairs $(x,K)$, where $K \in \mid
  \mathcal{P} \mid$ is an object of $\mathcal{P}$ and a morphism
  $(x,K)\to (y,L)$ is a pair $(\omega, i)$ where $i \in \mathcal{P}(K,L)$
  is a morphism in $\mathcal{P}$ from $K$ to $L$ and $\omega$ is a
  homotopy class of paths in $p^{-1}(L)$ from $y$ to $p^{-1}(i(x))$.
\end{dfn}
  \begin{dfn} \cite{AM,Att} The controlled homotopy groups $\pi_n^c(X,p)$ are defined to
  be the functor
  $$\pi_n^c(X):\mathcal{P}G_1(X)\to \mathcal{C}$$
  where $\mathcal{C}$ is tha category of pointed sets, groups or abelian
  groups defined by setting
  $$\pi_n^c(X,p)(x,K)=\pi_n(p^{-1}(K),x)$$
  and $\pi_n^c(\omega,i)$ is the composite of the change of basepoint
  isomorphism $\omega_*$ induced from $\omega$ and the homomorphism
  induced from the inclusion $i$.
\end{dfn}
\begin{dfn} \cite{AM,Att} The controlled homology $H^c_n(X,p)$ of a space controlled
  via $p:X \to Z$ is defined to be the pro-system $H_n(p^{-1}(B(r,z)))$
  via the maps $B(r,z) \to B(r+1,z)$.
  The controlled cohomology $H_c^n(X,p)$ of a space controlled via $p:X \to Z$
  is defined as a pro-system $H^n(p^{-1}(B(r,z))$ via the maps
  $B(r,z) \to B(r+1,z)$. 
\end{dfn}
\begin{dfn}\cite{AM,Att} If $(X,p)$ and $(Y,q)$ are spaces controlled over $Z$ by
  control maps $p$ and $q$ respectively, then $X$ is coextensive with $Y$
  if there exists an integer $m\ge 0$ so that if $p^{-1}(B_r(z))\ne \emptyset$
  then $q^{-1}(B_{r+m}(z))\ne \emptyset$ and the same with the roles of
  $p$ and $q$ reversed.
\end{dfn}
\begin{thm}[Whitehead Theorem] \cite{Att} If $f:(X,p) \to (Y,q)$ is a map of bounded
  simplicial complexes of bounded geometry, controlled by maps to Z
  of bounded geometry, then $f$ is a $bg$ homotopy equivalence if $(Y,q)$
  is coextensive with $(X,p)$ and for all $n\ge 0$, $f_*:\pi_n^c(X,p) \to
  f^!\pi_n^c(Y,q)$ is an isomorphism.
\end{thm}
\begin{thm}[Hurewicz Theorem] \cite{AM,Att} If $X$ is a simplicial complex
  of bounded geometry, then
  \par\noindent (i) $\pi_1^c(X)^{ab}=H_1^c(X)$
  \par\noindent (ii) If $\pi_i^c(X)=0$ for $i\le n-1$, $n \ge 2$ then,
  $H_i^c(X)=0$ for $i\le n-1$, and $H_n^c(X)=\pi_n^c(X).$
  \end{thm}
  \newtheorem{prop}{Proposition}[section]
  \begin{prop} We have an isomorphism of controlled homology and cohomology
    groups
  $$\tilde{H}_c^{p+k}(M^\nu)\overset{D}{\to}H^c_{n-p}(M)$$
  where the normal bundle $\nu$ and tangent bundle $\tau$ satisfy
  $\nu \oplus \tau=\epsilon^{n+k}$ where $n$ is the dimension of $M$. Here,
  $M^\nu$ is controlled over $M$ by the projection, and $M$ is controlled
  over itself by the identity.
\end{prop}
\begin{dfn} Let $E$ be a spectrum, $M$ a manifold of bounded geometry.
  The uniformly finite homology of $M$ with coefficients in $E$ is defined
  as
  $$E^{uf}_*(M)=[M^\nu:\Omega^{-*}E]^{Lip}$$
  the Lipschitz homotopy classes of Lipschitz maps from $M^\nu$ to
  $\Omega^{-*}E$ with uniform Lipschitz constant
  where $M^\nu$ is the S-dual of $M$ defined in the previous
  definition and $\Omega^{-1}X$ has the metric induced on the bar construction
  of $X$ by the product of the metrics on $X$. Again the same lack of
  independence of the homology on the choice of topological type of the
  spectrum holds as for Definition 2.28.
\end{dfn}

\begin{dfn} The surgery space \cite{Quinn} is defined as follows. Let $X$
  be a topological $n$-ad, $m \in \textbf{Z}$, and $w:\pi_1(X)\to \textbf{Z}_2$
  be a homomorphism. Define $\textbf{L}_m^k(X)$ or $\textbf{L}_m^s(X)$ as
  $\Delta$-set with $k$-simplices the set of topological surgery maps of
  $(k+n+3)$-ads of dimension $m+k$, $f:M \to K$ a (simple) homotopy equivalence on
  the $(k+3)rd$ face and with a reference map $h:K \to X$ with the orientation
  homomorphism $\pi_1(X)\to \textbf{Z}_2$ factoring through $w$ and
  $h(\partial_{k+4+j}K)\subseteq \partial_jX$. The first $k+2$ boundaries of
  this object are its boundaries as a $k$-simplex of the set. The $(k+3)$rd
  face, where it is a homotopy equivalence, is the connection with surgery,
  and the rest of the faces are there to map into the faces of $X$.
\end{dfn}
\begin{dfn} Define the complexity norm on the surgery space $\textbf{L}^s_m(X)$
  to be a norm assigning to each simplex the sum of the complexities of the
  surgery maps, being the number of simplices in the inverse image of each
  map, plus the number of simplices in a fixed triangulation of each space.
\end{dfn}
\begin{dfn} \cite{Ran2} The quadratic Poincar\'e n-ads over $\textbf{Z}[\pi]$ are the
  simplexes of -(k+1)-connected Kan complexes $\textbf{L}(\pi)$ $k\in
  \textbf{Z}$ such that $\Omega \textbf{L}_k=\textbf{L}_{k+1}$.
  Let $\sigma_*:G/TOP \to \textbf{L}_0$ the map associating to each singular
  simplex $\Delta \to G/TOP$ the quadratic Poincar\'e n-ad $\sigma_*(f,b)$
  over $\textbf{Z}$ of the normal map of manifold $n$-ads $(f,b)\to \Delta$
  that it classifies. This map induces the surgery obstruction isomorphism
  $$\sigma_*:\pi_*(G/TOP)\to \pi_*(\textbf{L}_0)$$
  We have $\textbf{L}[1/2]=BO[1/2]$ from Sullivan's characterization
  \cite{Su2} of stable
  topological bundles as $KO[1/2]$-oriented spherical fibrations.
\end{dfn}
\begin{dfn} Define $H^*_{uf}(X;\textbf{L})$ to be the group:
  $$[X:\textbf{L}(e)]^{bound.compl.}$$
  where $bound.compl.$ means homotopy classes of maps uniformly bounded in
  complexity norm with homotopies of bounded complexity norm.
\end{dfn}

\begin{prop}$H^*_{uf}(X;\textbf{L})\simeq [X:\textbf{L}(e)]^{Lip}$
\end{prop}
\textit{Proof} Take each surgery problem to its surgery obstruction. Then
Lipschitz homotopy classes of Lipschitz maps can be seen to be the same
as bounded complexity via the signature operator.
\begin{thm}[Narasimhan-Ramanan] \cite{NR} Let $G$ be a compact Lie group
and $n$ a positive integer. There exists a principal $G$-bundle $B$ and a
connection form $\gamma_1$ on $B$ so that for every principal $G$-bundle
$P$ with base of dimension $\le n$ and any connection form $\gamma$ on $P$,
one can find a bundle homomorphism $f$ of $P$ to $B$ so that the inverse
image of $\gamma_1$ by $f$ is $\gamma$.
\end{thm}
\newtheorem{lmm}{Lemma}[section]
\begin{lmm}[Gromov] \cite{Shubin} Let $X$ be a manifold of bounded geometry. There 
exists $\epsilon_0>0$ so that if $\epsilon \in (0,\epsilon_0)$ then there 
exists a countable covering of $X$ by balls of the radius $\epsilon$:
$X=\cup B(x_i, \epsilon)$ so that the covering of $X$ by the balls
$B(x_i, 2\epsilon)$ with the double radius and the same centers has 
a finite maximal number of balls with nonempty intersection in this covering.
\end{lmm}
\begin{lmm}\cite{Engel} Let a covering $\{U_\alpha\}$ of a manifold $M$ of 
bounded geometry  be given having finite multiplicity. Then there exists
a coloring of the subsets of $U_\alpha$  with finitely many colors such
that no two intersecting subsets have the same color.
\end{lmm}
\begin{thm} Let $X$ be a manifold of bounded geometry, and $(E,\nabla)$
a vector bundle of bounded geometry over $X$. Then there is a Lipschitz
map $f:X \to BO(n)$ classifying $E$ 
so that if $\gamma$ is the universal connection, $f^*(\gamma)=\nabla$.
\end{thm}
\textit{Proof:} Since $X$ and $E$ have bounded geometry, we can find a
uniformly locally finite cover of $X$ by balls of a fixed radius together
with a subordinate partition of unity whose derivatives are all uniformly
bounded and such that over each coordinate ball $E$ is trivalized by
a synchronous framing. 

Now we  use the lemma to color the coordinate balls with finitely many
colors so that no two balls with the same color intersect. This gives
a partition of the coordinate balls into $N$ families $U_1,...,U_N$
so that any $U_i$ is a collection of disjoint balls, and we get a 
corresponding subordinate partition of unity $1=\phi_1+...+\phi_N$
with uniformly bounded derivatives. Since $E$ is trivial over each 
$U_i$, we can apply the Narasimhan-Ramanan theorem to $\nabla\mid U_i$
to get a classifying map $f_i:U_i \to BO$ so that $f^*(\gamma)=\nabla \mid
U_i$ and piece the maps together with the partition of unity.
\begin{thm} Let $X$ be a manifold of bounded geometry. Let
$KO_{uf}(X)$ be the uniformly finite KO-theory of $X$. Then 
$$KO_{uf}(X)=[X:BO]^{Lip}$$
where $[X:Y]^{Lip}$ is the group of Lipschitz homotopy classes of Lipschitz maps with uniform Lipschitz constant,
and $BO$ is equipped with the universal connection. Here $E=BO$ is the spectrum
with $E_0=\bigcup_{n=0}^\infty BO(n)$ so that it is a CW complex with finite
$n$-skeleta.
\end{thm}
\begin{rk} There is another metric on $BO$, associated to $O(n)$ where $n$
  gets larger and one uses the operator norm. In this metric, all elements of
  $KO(S^n)$ have bounded norm. In the metric of Theorem 2.11, only finitely
  many do.
  \end{rk}
\textit{Proof:} Let $(\xi, \theta)$ denote an element of $KO_{uf}(X)$. 
There is a classifying map $g:X \to BO$, which maps to a finite-dimensional 
skeleton of $BO$ classifying $\xi$ and by the 
Narasimhan-Ramanan theorem there is a map $g^\prime$ homotopic to $g$
with $g^{\prime*}(\gamma_1)=\theta$, where $\gamma_1$ is the universal 
connection on $BO$. If we bound the derivatives of $g^\prime$, we
bound the Christoffel symbols of $\theta$. Then
clearly a Lipschitz map pulls back the universal connection to $\theta$,
so that the connection form is in $L^\infty$. Addition is direct sum, and 
inverses exist, since if we take the ordinary inverse of $\xi$ in 
$KO_{uf}(X)$ the direct sum with $\xi^{-1}$ has a globally flat connection
when it is embedded in a high dimensional Euclidean space \cite{CS}.
\begin{dfn} \cite{Moerd}. Let $(M, \mathcal{F})$ be a foliated manifold. The holonomy groupoid
$H=Hol(M,\mathcal{F})$ is a smooth groupoid with $H_0=M$ as the space
of objects. If $x,y \in M$ are two points on different leaves there are
no arrows between $x$ and $y$ in $H$. If $x$ and $y$ lie on the same leaf
$L$, an arrow $h:x \to y$ in $H$ is an equivalence class $h=[\alpha]$ of 
smooth paths $\alpha:[0,1] \to L$ with $\alpha(0)=x$ and $\alpha(1)=y$. 
To explain the equivalence relation, let $T_x$ and $T_y$ be small q-disks 
through $x$ and $y$ transverse to the
leaves of the foliation. If $x^\prime \in T_x$ is a point sufficiently close
to $x$ on a leaf $L^\prime$, then $\alpha$ can be ``copied'' inside
$L^\prime$ to give a path $\alpha^\prime$ near $\alpha$ with endpoint
$y^\prime \in T_y$. In this way one obtains the germ of a diffeomorphism
from $T_x$ to $T_y$ sending $x$ to $y$ and $x^\prime$ to $y^\prime$. This
germ is called the holonomy of $\alpha$ and denoted $hol(\alpha)$. Two
paths $\alpha$ and $\beta$ from $x$ to $y$ in $L$ are equivalent, i.e.
define the same arrow $x \to y$,  if and
only if $hol(\alpha)=hol(\beta)$. We obtain a well defined smooth groupoid 
$H=Hol(M,\mathcal{F})$. This groupoid is a foliation groupoid, and the 
(discrete) isotropy group $H_x$ at $x$ is called the holonomy group of the 
leaf through $x$. 
\end{dfn}
Alternatively \cite{Connes}, let $(V,F)$ be a foliated manifold of codimension
$q$. Every point $x \in V$ has a neighborhood $U$ and a system of local
coordinates $(x^j)_{j=1,...,dim(V)}$, which is called a foliation chart, so
that the partition of $U$ into connected components of leaves, called
plaques, corresponds to the partition of $\textbf{R}^{dim(V)}=
\textbf{R}^{dim(F)}\times \textbf{R}^{codim(F)}$ into the parallel affine
subspaces $\textbf{R}^{dim(F)} \times pt$. Given any $x \in V$ and a small enough open set $W \subset V$ containing
$x$, the restriction of the foliation $F$ to $W$ has, as its leaf space, an
open set of $\textbf{R}^q$, which we shall call for short a transverse 
neighborhood of $x$. In other words this open set $W/F$ is the set of plaques
around $x$. Now, given a leaf $L$ of $(V,F)$ and two points $x,y \in L$ of 
this leaf, any simple path $\gamma$ from $x$ to $y$ on the leaf $L$ 
uniquely determines a germ $h(\gamma)$ of a diffeomorphism from a transverse
neighborhood of $x$ to a transverse neighborhood of $y$. One can obtain 
$h(\gamma)$, for instance, by restricting the foliation $F$ to a neighborhood
$N$ of $\gamma$ in $V$ sufficiently small to be a transverse neighborhood
of both $x$ and $y$ as well as of any $\gamma (t)$. The germ of diffeomorphism
$h(\gamma)$ thus obtained only depends upon the homotopy class $\gamma$
only depends upon the homotopy class of $\gamma$ in the fundamental groupoid
of the leaf $L$, and is called the holonomy of the path $\gamma$. The 
holonomy groupoid of a leaf $L$ is the quotient of its fundamental groupoid
by the equivalence relation which identifies two paths $\gamma$ and 
$\gamma^\prime$ from $x$ to $y$ in $L$ if and only if $h(\gamma)=h(\gamma^\prime)$. 
The holonomy covering $\tilde{L}$ of a leaf is the covering of $L$ associated
to the normal subgroup of its fundamental group $\pi_1(L)$ given by paths 
with trivial holonomy. The holonomy groupoid of the foliation is the union
$G$ of the holonomy groupoids of its leaves. Given an element $\gamma$ of $G$,
we denote by $x=s(\gamma)$ the origin of the path $\gamma$, by $y=r(\gamma)$
its end point, and $r$ and $s$ are the range and source maps.
\begin{dfn}\cite{Moerd}
For a smooth groupoid $G$, the nerve of $G$ is the simplicial set whose 
$n$-simplices are strings of composable arrows in $G$. This set is denoted
$G_n$. The space $G_n$ is the fibered product $G_1 \times_{G_0} ... \times_{G_0} G_1$, hence has the natural structure of a smooth simplicial manifold, that
is a simplicial set which has the structure of a smooth manifold.
 Thus $G_\cdot$ is a simplicial manifold whose geometric realization is 
denoted $BG$ and is called the classifying space of $G$.
\end{dfn}
The following alternative definition of $BG$ is from \cite{Connes}.
The space $BG$ is only defined up to homotopy as the quotient space of a free
and proper action of $G$ on contractible spaces. More precisely, a right
action of $G$ on topological spaces is given by a topological space $Y$,
a continuous map $s_Y:Y \to G^{(0)}=M$ (to the zero dimensional simplices)
 and a continuous map
$Y \times_s G \to^\circ Y$ where $Y \times_s G=\{(y,\gamma) \in Y \times G;
s_Y(y)=r(\gamma)\}$ such that $(y \circ \gamma_1) \circ \gamma_2=
y\circ(\gamma_1 \gamma_2)$ for any $y \in Y$, $\gamma_1 \in G, \gamma_2 \in G$
with $r(\gamma_1)=s_Y(y)$, $s(\gamma_1)=r(\gamma_2).$

In other words the map which to each $t \in G^{(0)}$ assigns the topological
space $Y_t=s_Y^{-1}\{t\}$ and to each $\gamma \in G$, $\gamma: t \to t^\prime$,
assigns the diffeomorphism $y \to y \circ \gamma$ from $Y_{t^\prime}$ to 
$Y_t$, is a contravariant functor from the small category $G$ to the category
of topological spaces. We shall say that such an action of $G$ on $Y$ is 
free iff for any $y \in Y$ the map $\gamma \in G^{s_Y(y)} \to y \circ \gamma
\in Y$ is injective, and that it is proper iff the map $(y,\gamma) \to
(y, y\circ \gamma)$ is proper. A free and proper action of $G$ on $Y$ is 
the same thing as a principal $G$-bundle on the quotient space 
$Z=Y/G$, the quotient of $Y$ by the equivalence relation $y_1 \sim y_2$
iff there exists $\gamma \in G$ so that $y_1 \circ \gamma = y_2$. We say that
the action of $G$ on $Y$ has contractible fibers iff the fibers of 
$s_Y: Y \to G^{(0)}$ are contractible. Then as for groups, the classifying
space $BG$ is the quotient $EG/G$ of a principal $G$-bundle $Y=EG$ with 
contractible fibers. 
\begin{dfn} \cite{BCH} 
Let $X$ be a simplicial complex. A cosheaf of abelian groups $\mathcal{A}$
on $X$ consists of the following: For each simplex on $X$ an abelian
group $A_\sigma$. For each inclusion of simplices $\eta \subseteq \sigma$
a homomorphism of abelian groups $\phi_\eta^\sigma:A_\sigma \to A_\eta$ with
$\phi_\tau^\sigma = \phi_\tau^\eta \phi_\eta^\sigma$ whenever 
$\tau \subseteq \eta \subseteq \sigma$, and with 
$\phi^\sigma_\sigma=id$ for each $\sigma$.
\end{dfn}
\begin{dfn}\cite{BCH} If $\eta$ and $\sigma$ are oriented simplices and 
if $\eta$ is a codimension 1 face of $\sigma$ define an incidence number
$$[\eta:\sigma]=\begin{cases} +1 \mbox{ if the orientation on }\eta \mbox{  is inherited from the one from } \sigma\\-1 \mbox{ otherwise } \end{cases}$$
\end{dfn}
\begin{dfn}\cite{BCH} Let $X$ be a simplicial complex, $\mathcal{A}$ a cosheaf
on $X$. Denote by $C_n(X;\mathcal{A})$ the abelian group whose 
elements are all formal sums
$$\sum_{dim \sigma=n} a_\sigma[\sigma]$$
where $\sigma$ ranges over all simplices of dimension $n$ and 
$a_\sigma \in A_\sigma$. Define the boundary
$$\partial:C_{n+1}(X;\mathcal{A}) \to C_n(X;\mathcal{A})$$
by orienting $X$ and using the formula
$$\partial(a_\sigma[\sigma])=\sum_{dim(\eta)=dim(\sigma)-1}
[\eta : \sigma]\phi^{\sigma}_{\eta}(a_\sigma)[\eta]$$
The homology groups of this complex are the homology groups of $X$ with
coefficients in the cosheaf $\mathcal{A}$.
\end{dfn}
\begin{rk} Because we are assuming that the foliation $\mathcal{F}$ is a 
foliation of a compact manifold $M$, we do not need norms on the groups of 
the cosheaf in the above definition. 
\end{rk}
\section{Integrability}
Given $\pi:E \to X$ of a vector bundle of dimension $n$, the paracompactness 
of $X$ makes it possible to find a countable open cover $\{U_i\}_{i=1}^\infty$
of $X$ and a subordinate partition of unity $\{\lambda_i\}_{i=1}^\infty$ and
a continuous $\phi_i:E \mid U_i \to \textbf{R}^n$ mapping each fiber 
isomorphically onto $\textbf{R}^n$. Express a separable real Hilbert space
 $\mathcal{H}$ as a countable orthogonal direct sum of copies of 
$\textbf{R}^n$ and let $\psi_i:\textbf{R}^n \to \mathcal{H}$ be the inclusion
of $i^{th}$ summand. Then $\phi:E \to \mathcal{H}$ defined by
$$\phi(e)=\sum_{i=1}^\infty \lambda_i(\pi(e))\cdot \psi_i \circ \phi_i(e)$$ 
is continuous and sends each fiber $\pi^{-1}(x)$ linearly and one-one onto
an $n$-dimensional subspace of $\mathcal{H}$. Thus $g_E(x)=\phi(\pi^{-1}(x))$
defines a continuous map $g_E:X \to BGL_n(\textbf{R})$ \cite{BC}.
\begin{dfn}\cite{BC} Let $E$ be an $n$-dimensional vector bundle over $X$ and let 
$g_E:X \to BGL_n(\textbf{R})$ be the map classifying $E$. Then
the graded subring
$$Pont^*(E)=g_E^*(H^*(BGL_n(\textbf{R});\textbf{R})) \subset H^*(X;\textbf{R})$$
is called the Pontrjagin ring of $E$.
\end{dfn}
\begin{dfn} Let $(E,\nabla)$ be an $n$-dimensional vector bundle of bounded geometry
over a smooth manifold $M$ of bounded geometry and let 
$$g^{(\infty)}_{(E,\nabla)}:M \to BO_n$$
be the classifying map for $(E,\nabla)$. Then the $L^\infty$ Pontrjagin
ring of $E$ is the graded ring
$$Pont^{(\infty),*}(E)=g_{(E,\nabla)}^{(\infty)*}(H^*(BO_n;\textbf{R})) \subset H_\beta^*(M; \textbf{R})\simeq H_{uff}^*(M;\textbf{R})$$
generated by the $L^\infty$ Pontrjagin classes of E, 
where the isomorphism between bounded de Rham theory and uniformly finite 
cohomology is given in \cite{AttB}.
\end{dfn}
\begin{dfn} \cite{Ran} Let $C$ be a chain complex, the suspension of $C$,
$SC$ is defined by
$$d_{SC}: SC_r=C_{r-1} \to SC_{r-1}=C_{r-2}$$
The $n$-fold suspension $S^nC$ is defined to be
$$(S^nC)_r=C_{r-n}$$
\end{dfn}
\begin{dfn} \cite{Ran} Let $\mathcal{A}$ be an additive category.
We denote by $\mathcal{B}(\mathcal{A})$ the 
additive category of finite chain complexes in $\mathcal{A}$ and chain maps.
\end{dfn}
\begin{dfn} \cite{Ran} (i) Let $K$ be a simplicial complex. 
 An object $M$ in an additive category $\mathcal{A}$ is 
$K$-based if it is expressed as the direct sum
$$M=\sum_{\sigma \in K} M(\sigma)$$
of objects $M(\sigma)$ in $\mathcal{A}$, corresponding to simplices $\sigma$
of $K$ so that $\{\sigma \in K \mid
M(\sigma) \ne 0\}$ is finite. A morphism $f:M \to N$ of K-based objects is 
a collection of morphisms in $\mathcal{A}$
$$f==\{f(\tau,\sigma):M(\sigma) \to N(\tau) \mid \sigma,\tau \in K\}$$
\par (ii) Let $\mathcal{A}_*(K)$ 
be the additive category of K-based objects $M$ in $\mathcal{A}$ with
morphisms $f:M \to N$ such that $f(\tau,\sigma):M(\sigma) \to N(\tau)$
is 0 unless $\tau \ge \sigma$ so that
$$ f(M(\sigma))\subseteq\sum_{\tau \ge \sigma} N(\tau)$$
\par (iii) Let $\mathcal{A}^*(K)$ be the additive category of K-based objects
$M$ in $\mathcal{A}$ with morphisms $f:M \to N$ such that 
$f(\tau,\sigma):M(\sigma)\to N(\tau)$ is 0 unless $\tau \le \sigma$ so that
$$f(M(\sigma))\subseteq \sum_{\tau\ge\sigma}N(\tau)$$
\par (iv) Forgetting the K-based structure defines the covariant assembly
functor:
$$\mathcal{A}_*(K) \to \mathcal{A}; M \to M_*(K)=\sum_{\sigma \in K} M(\sigma)$$
\end{dfn}
\begin{dfn}\cite{Ran}(i) Let $\mathcal{A}_*[K]$ be the 
additive category with objects the covariant additive functors
$$M:K \to \mathcal{A};\sigma \to M(\sigma)$$
such that $\{\sigma \in K \mid M[\sigma] \ne 0\}$ is finite. The morphisms
are the natural transformations of such functors.
\par\noindent (ii) Let $\mathcal{A}^*[K]$ be the additive category with
objects the contravariant additive functors
$$M:K\to \mathcal{A}; \sigma \to M[\sigma]$$
such that $\{\sigma \in K \mid M[\sigma] \ne 0\}$ is finite. The morphisms
are the natural transformations of such functors.
\end{dfn}
\begin{dfn}\cite{Ran} Define a contravariant functor 
$T:\mathcal{A}^*[K]\to\mathcal{B}(\mathcal{A}^*(K))$ by sending an object
$M$ to the chain complex $TM$ with
$$(TM)_r(\sigma)=T(M[\sigma])_{r-\mid\sigma\mid}$$
$$d_{TM}(\tau,\sigma)=(-1)^iT(M[\tau] \to M[\sigma]):(TM)_r(\sigma)\to
(TM)_{r-1}(\tau)$$
if $\sigma \ge \tau$, $\mid \sigma \mid=\tau+1$, $\tau=\partial_i \sigma$.
\end{dfn}
\begin{dfn}\cite{Ran} Define the covariant assembly functor for a simplicial 
complex $K$
$$\mathcal{B}(\mathcal{A})_*[K]=\mathcal{B}(\mathcal{A}_*[K]) \to \mathcal{B}(\mathcal{A})^*(K)=\mathcal{B}(\mathcal{A}^*(K));C \to C_*[K]$$
by sending a finite chain complex $C$ in  $\mathcal{A}_*[K]$
 to the finite chain complex $C_*[K]$ in  $\mathcal{A}_*(K)$ with
$$C_*[K]_r=\sum_{\sigma \in K} C[\sigma]_{r-\mid \sigma \mid}, C_*[K](\sigma)=
S^{-\mid \sigma \mid} C[\sigma]$$
$\sigma \in K$. The assembly is the total complex of the double complex in
$\mathcal{A}$ defined by
$$C_*[K]_{p,q}=\sum_{\sigma \in K, \mid \sigma \mid=-p} C[\sigma]_q$$
$$d^\prime:C_*[K]_{p,q}\to C_*[K]_{p-1,q}; c[\sigma]\to \sum_i (-1)^i\partial_i c[\sigma]$$
$$d^{''}:C_*[K]_{p,q} \to C_*[K]_{p,q-1}; c[\sigma] \to d_{C[\sigma]}(c[\sigma])$$
the sum of $d^\prime$ being taken over all the elements $\partial_i \sigma \in K_*(\sigma)$ with $\partial_i: C[\sigma]\to C[\partial_i \sigma]$
 the chain map induced by the inclusion $\partial_i\sigma \to \sigma$.
\end{dfn}
\begin{dfn}\cite{Ran} A chain duality $(T,e)$ on an additive category $\mathcal{A}$ is
a contravariant additive functor $T:\mathcal{A} \to \mathcal{B}(\mathcal{A})$ 
together with a natural transformation
$$e:T^2 \to 1:\mathcal{A}\to \mathcal{B}(\mathcal{A})$$
such that for each object $A$ in $\mathcal{A}$
\begin{enumerate}
\item $e(T(A))\cdot T(e(A))=1:T(A)\to T^3(A)\to T(A)$.
\item $e(A):T^2(A)\to A$ is a chain equivalence.
\end{enumerate}
\end{dfn}
\begin{dfn}\cite{Ran} Use the standard free $\textbf{Z}[\textbf{Z}_2]$-module 
resolution of $\textbf{Z}$ to define for any finite chain complex $C$
in $\mathcal{A}$ the $\textbf{Z}$-module chain complex
$$W_{\%}C=W\otimes_{\textbf{Z}[\textbf{Z}_2]}(C\otimes_\mathcal{A} C)$$
\end{dfn}
\begin{dfn}\cite{Ran} An $n$-dimensional quadratic Poincar\'e complex in $\mathcal{A}$
$(C,\psi)$ is a finite chain complex $C$ in $\mathcal{A}$ together with an
$n$-cycle $\psi \in (W_\%)_n$ such that the chain map
$$(1+T)\psi_0:C^{n-*}\to C$$
is a chain equivalence in $\mathcal A$.
\end{dfn}
\begin{dfn}\cite{Ran} Let $f:C \to D$ be a chain map. The algebraic mapping cone of $f$,
$C(f)$ is defined by 
$$d_{C(f)}=\begin{pmatrix} d_D & (-1)^{r-1}f\\ 0 & d_C \end{pmatrix}$$
$$C(f)_r=D_r\oplus C_{r-1} \to C(f)_{r-1}$$
\end{dfn}
\begin{dfn}\cite{Ran} A chain map $f:C \to D$ of finite chain complexes in $\mathcal{A}$
induces $\textbf{Z}[\textbf{Z}_2]$-module chain map
$$f \otimes f:C\otimes_{\mathcal{A}}C \to D\otimes_{\mathcal{A}} D$$
and hence a $\textbf{Z}$-module chain map
$$f_{\%}: W_\%C \to W_\%D$$
\end{dfn}
\begin{dfn}\cite{Ran} An $n+1$-dimensional quadratic Poincar\'e pair in $\mathcal{A}$
$$(f:C\to D, (\delta\psi,\psi))$$
is a chain map $f:C \to D$ of finite chain complexes together with an 
$n+1$-cycle $(\delta\psi,\psi)\in C(f_\%)_{n+1}$ such that the chain map
$$(1+T)(\delta\psi_0, \psi_0):D^{n+1-*}\to C(f)$$
is a chain equivalence.
\end{dfn}
\begin{dfn}\cite{Ran} A cobordism of $n$-dimensional quadratic Poincar\'e complexes
$(C,\psi)$, $(C^\prime,\psi^\prime)$ is an $n+1$-dimensional Poincar\'e pair
$$(f,f^\prime):C\oplus C^\prime \to D, (\delta\psi,\psi\oplus -\psi^\prime)$$
\end{dfn}
\begin{dfn}\cite{Ran} A subcategory $\mathcal{C} \subseteq \mathcal{B}(\mathcal{A})$
is closed if it is a full additive category so that the algebraic mapping
cone $C(f)$ of any chain map $f:C \to D$ in $\mathcal{C}$ is an object
of $\mathcal{C}$. A chain complex $C$ in $\mathcal{A}$ is $\mathcal{C}$-contractible if it belongs to $\mathcal{C}$. An $n$-dimensional quadratic complex
$(C,\psi)$ in $\mathcal{A}$ is $\mathcal{C}$-contractible if $C$ and 
$C^{n-*}$ are $\mathcal{C}$-contractible. An $n$-dimensional quadratic 
complex $(C,\psi)$ in $\mathcal{A}$ is $\mathcal{C}$-Poincar\'e if the chain
complex
$$\partial C=S^{-1}C((1+T)\psi_0:C^{n-*} \to C)$$
is $\mathcal{C}$-contractible.
\end{dfn}
\begin{dfn}\cite{Ran} Let $\Lambda=(\mathcal{A}, \mathcal{B},\mathcal{C})$ be a triple
of additive categories, where $\mathcal{A}$ has chain duality
$T:\mathcal{A} \to \mathcal{B}(\mathcal{A})$ and a pair 
$(\mathcal{B}, \mathcal{C} \subseteq \mathcal{B})$ of closed subcategories
of $\mathcal{B}(\mathcal{A})$ so that for any object $B$ of $\mathcal{B}$
\begin{enumerate}
\item The algebraic mapping cone $C(1: B \to B)$ is an object of $\mathcal{C}$.
\item The chain equivalence $e(B):T^2(B) \to B$ is a $\mathcal{C}$-equivalence.
\end{enumerate}
Then $\Lambda$ is said to be an algebraic bordism category.
\end{dfn}
\begin{dfn}\cite{Ran} For any additive category with chain duality $\mathcal{A}$ there 
is defined an algebraic bordism category
$$\Lambda(\mathcal{A})=(\mathcal{A}, \mathcal{B}(\mathcal{A}), \mathcal{C}(\mathcal{A}))$$
with $\mathcal{B}(\mathcal{A})$ the category of finite chain complexes 
in $\mathcal{A}$ and $\mathcal{C}(\mathcal{A})\subseteq \mathcal{B}(\mathcal{A})$ the subcategory of contractible chain complexes.
\end{dfn}
\begin{dfn}\cite{Ran} Let $\Lambda=(\mathcal{A}, \mathcal{B}, \mathcal{C})$ be an 
algebraic bordism category. An $n$-dimensional quadratic complex $(C, \psi)$
in $\mathcal{A}$ which is $\mathcal{B}$-contractible and $\mathcal{C}$-Poincar\'e. The quadratic $L$-group $L_n(\Lambda)$ is the cobordism group of 
$n$-dimensional quadratic complexes in $\Lambda$.
\end{dfn}
\begin{dfn} Let $K$ be a simplicial complex of bounded geometry.
Let  $\mathcal{A}^{uf}(K)$ be the additive category of 
$K$-based objects in $\mathcal{A}$ which fall into a finite number
of isomorphism types in every ball of fixed radius $R$, depending on the 
radius, in $K$.
\end{dfn}
\begin{dfn} Let $\mathcal{A}(R)^{uf}_*(K)$ be the additive category of 
$K$-based objects in $\mathcal(R)$ which fall into a finite number
of types in a ball of fixed radius with morphisms $f:M \to N$ so 
that $f(\tau, \sigma)=0:M(\sigma) \to N(\tau)$ unless $\tau \ge \sigma$
so that $f(M(\sigma))\subseteq \sum_{\tau \ge \sigma} N(\tau)$.
\end{dfn}
\begin{dfn} Define the algebraic bordism category of local, uniformly finite,
finitely generated free $(R,K)$-modules, where $R$ is a ring and $K$ a
simplicial complex of bounded geometry as
$$\Lambda(R)_*^{uf}(K)=(\mathcal{A}^{uf}(R,K), \mathcal{B}^{uf}(R,K), \mathcal{C}(R)_*^{uf}(K))$$
where $\mathcal{B}^{uf}(R,K)$ is the category of finite chain complexes of 
f.g. free uniformly finite $(R,K)$-modules. An object $\mathcal{C}(R)_*^{uf}(K)$
is a finite f.g. free uniformly finite $(R,K)$-module chain complex $C$ 
such that each $[C][\sigma]$, $\sigma \in K$, is a contractible f.g. free
$R$-module chain complex.
\end{dfn}
\begin{dfn} Let $M$ be a manifold of bounded geometry. We define an object
which we will refer to as the algebraic uniformly finite homology of $M$ with 
coefficients in the $L$-spectrum
$$H_n^{uf}(M;\textbf{L})=L_n(\Lambda(\textbf{Z})_*^{uf}(M))$$
\end{dfn}
\begin{dfn} Define the complex algebraic uniformly finite homology of $M$ with
coefficients in the $L$-spectrum as
$$H_n^{uf}(M;\textbf{L}(\textbf{C}))=L_n(\Lambda(\textbf{C})_*^{uf}(M))$$
\end{dfn}
\begin{prop}[Ranicki] \cite{Ran} An algebraic bordism category $\Lambda=(\mathcal{A},
  \mathcal{B},\mathcal{C})$ and a locally finite ordered simplicial
  complex $K$ determine an algebraic bordism category
  $$\Lambda^*(K)=(\mathcal{A}^*(K), \mathcal{B}^*(K), \mathcal{C}^*(K))$$
  $$\Lambda_*(K)=(\mathcal{A}_*(K), \mathcal{B}_*(K), \mathcal{C}_*(K))$$
\end{prop}
\begin{dfn}\cite{Ran} Let $K$ and $L$ be $\Delta$-sets. Then the geometric product of
  $K$ and $L$, $K \otimes L$ is the $\Delta$-set with one $p$-simplex
  for each equivalence class of triples
  $$(\mbox{m-simplex }\sigma\in K,\mbox{ n-simplex }\tau\in L, \mbox{ p-simplex
  }\rho\in\Delta^m\otimes\Delta^n)$$
  subject to the equivalence relation generated by
  $$(\sigma, \tau,\rho) \sim (\sigma^\prime, \tau^\prime, \rho^\prime)
  \mbox{ if there exist }\Delta\mbox{-maps }f:\Delta^m\to\Delta^{m^\prime}$$
  $$g:\Delta^n\to\Delta^{n^\prime}\mbox{ such that }\sigma=f^*\sigma^\prime,
  \tau=g^*\tau^\prime, (f\otimes g)_*(\rho)=(\rho^\prime)$$
  \end{dfn}
\begin{dfn} Let $K$ be a $\Delta$-set, $L$ be a $\Delta$-set with a
  complexity function assigned to its simplices, let the $\Delta$-set
  $L^{K,bound. compl.}$ be the $\Delta$-set with
  $(L^{K,bound. compl.})^{(n)}$ the set of $\Delta$-maps of bounded complexity
  $K\otimes\Delta^n \to L$, with $\partial_i$ induced from $\partial_i:
  \Delta^{n-1} \to \Delta^n$.
  \end{dfn}
  \begin{prop} Let $M$ be a manifold of bounded geometry. Then
  $$L_n(\Lambda(\textbf{Z})_*^{uf}(M))=[M^\nu:\Omega^{-n}\textbf{L}(e)]^{Lip}$$
\end{prop}
\textit{Proof}; In section 13 of \cite{Ran} we have for any additive category
$\mathcal{A}$ with duality with algebraic bordism category
$$\Lambda(\mathcal{A})=(\mathcal{A}, \mathcal{B}(\mathcal{A}),\mathcal{C}(\mathcal{A}))$$
Let $\textbf{L}_n(\Lambda)$ be the $\Delta$-set with $m$-simplexes the
$n$-dimensional quadratic complexes in $\Lambda^*(\Delta^m)$, with the zero
complex as $m$-simplex $\emptyset$. The face maps are indeced from the
standard embeddings $\partial_i:\Delta^{m-1} \to \Delta^m$ via the functors
$$(\partial_i)^*:\Lambda^*(\Delta^m)\to \Lambda^*(\Delta^{m-1})$$
We then imitate the proof of 13.7 in \cite{Ran}:
$$\textbf{L}_n(\Lambda^*_{uf}(M))^{(p)}=\{\mbox{n-dimensional quadratic
  complexes in }\Lambda^*_{uf}(M)^*(\Delta^p)\}$$
$$(\textbf{L}_\cdot(\Lambda)^{M_+, bound. compl})^{(p)}=\{\mbox{n-dimensional
  quadratic complexes in }\Lambda^*_{uf}(M \otimes \Delta^p)\}$$
$$=\{\mbox{Bounded complexity maps }(M \otimes \Delta^p)_+ \to \textbf{L}_n
(\Lambda)\}$$
by definition. 
We get a representation in cohomology
$$\textbf{L}_n(\Lambda_{uf}^*(M))=M_+^{\textbf{L}_n(\Lambda),bound. compl.}=
\textbf{H}^n_{uf}(M;\textbf{L}_\cdot(\Lambda))$$
  We then take the Spanier-Whitehead dual
  $$\textbf{L}(\Lambda_*^{uf}(M))^{(p)}=\{\mbox{n-dimensional quadratic
    complexes in }\Lambda_*^{uf}(M)^*(\Delta^p)\}$$
  We obtain
  $$\textbf{H}_n^{uf}(M;\textbf{L}_\cdot(\Lambda))^{(p)}=\{\mbox{(n-m)-dimensional
    quadratic complexes in }\Lambda^*_{uf}(M^\nu \otimes \Delta^p)\}$$
  $$=\textbf{H}^{n-m}_{uf}(M^\nu;\textbf{L}_\cdot(\Lambda))^{(p)}=
  \textbf{L}_{n-m}(\Lambda^*_{uf}(M^\nu))^{(p)}$$
  $$=\{\mbox{Bounded complexity maps }M^\nu \otimes \Delta^p \to
  \textbf{L}_{n-m}(\Lambda)\}$$
  
\begin{thm} Let $X$ be a $bg$ manifold. Then
  $$H^{uf}_*(X;\textbf{L})_{(p)}=KO^{uf}_*(X)_{(p)}$$
  where $p$ is an odd prime. 
\end{thm}
\textit{Proof} Consider $H^*_{uf}(X;\textbf{L})=[X:\textbf{L}(e)]^{bound.compl.}$.
Now $KO^*_{uf}(X)=[X:BO]^{Lip}$, and we have $\textbf{L}(e)_{(p)}=BO_{(p)}$.
The bounded complexity norm and the universal connection on $BO$ are
compatible, so we have
$$H^*_{uf}(X;\textbf{L})_{(p)}\simeq KO^*_{uf}(X)_{(p)}$$
Now we have Poincar\'e duality for $KO^*_{uf}(X)$ and $H^*_{uf}(X;\textbf{L})$
and we have the result.

One of the uses of the above Theorem is to relate L-theory
to KO-theory, although there are other ways to achieve this \cite{Ran},
Chapter 16.
\newtheorem{cor}{Corollary}[section]
\begin{cor}
$$H_n^{uf}(M;\textbf{L})\otimes \textbf{Q}=KO_*^{uf}(X)\otimes \textbf{Q}$$
\end{cor}
\textit{Proof} This follows directly from the theorem,
and the fact that in the proof of that theorem, the surgery spectrum is 
KO-theory away from the prime 2.
\begin{cor}
$$H_n^{uf}(M;\textbf{L}(\textbf{C}))=K_*^{uf}(X) \otimes \textbf{C}$$
\end{cor}
\newtheorem{conj}{Conjecture}[section]
\begin{conj}[Novikov Higher Signature Conjecture]
The assembly map
$$H_*(B\pi;\textbf{L})\otimes \textbf{Q}\to L_*(\pi)\otimes \textbf{Q}$$
is injective.
This is the same as the following: If $M$ is a manifold of dimension
$n$ and $\pi$ is a group so that there is a map $i:\pi_1(M) \to \pi$
and hence a map $F:M \to B\pi$ then if $\gamma \in H^{n-4i}(B\pi;\textbf{Q})$
then
$$\langle f^*(\gamma) \cup L_i(M),[M]\rangle\in\textbf{Q}$$
is homotopy invariant.
\end{conj}
\begin{dfn}\cite{Dran} A metric space $X$ is uniformly contractible if 
  there exists a function $S(r)>0$ such that every ball $B_r(x)$ of radius $r$
  around any point $x$ is contractible to a point in $B_{S(r)}(x)$, the ball of
  radius $S(r)$ around $x$.
  \end{dfn}
\medskip

Although L-theory and KO-theory are the same tensored with $\textbf{Z}[1/2]$
we will need to tensor with $\textbf{Q}$ to use the Pontrjagin character.
In addition, Bott and Heitsch \cite{BottH} have shown that the 
Integrability Theorem fails for integral Pontrjagin classes and also 
for $p$-torsion where $p$ is odd, so it cannot be true for characteristic 
classes with coefficients in $\textbf{Z}[1/2]$.

Let $E$ be a subbundle of the tangent bundle of a $C^\infty$ compact manifold $M$. Then we will 
show that if $E$ is integrable the Pontrjagin ring $Pont^*(TM/E)$ vanishes 
above dimension $2k$ where $k$ is the dimension of $TM/E$ \cite{Bott}. Suppose
$E$ is integrable.
Let $\mathcal{F}$ be the foliation on a manifold $M$ corresponding to $E$.
Then the tangent bundle of $\mathcal{F}$ is $E$ and the normal bundle
$\nu(\mathcal{F})$ is the set of vectors perpendicular to $E=T(\mathcal{F})$, 
or $TM/E$. 
Consider the smooth groupoid $G$ of $\mathcal{F}$
 and over each point in the groupoid $x=(x_1,y_1)$ take the corresponding leaf
$S_x$ through $(x_1,y_1)$. 
\begin{thm}[Bott Integrability]\cite{Bott}
Let $Pont^*(TM/E)$ be the Pontrjagin ring of $TM/E$. Then if E is integrable,
for $q>2k$, $Pont^q(TM/E)=0$.
\end{thm}
\textit{Proof:}
Let $BG$ be the classifying space of the foliation groupoid $G$ of 
$\mathcal{F}$ and consider the homology of $BG$ with coefficients
in the cosheaf $\mathcal{L}^{uf}(S_x)$ of L-homology groups of the leaves of
$\mathcal{F}$:
$H_*(BG; \mathcal{L}^{uf}(S_x))$. This group is the set of the normal
invariants of the ``blocked surgery'' of \cite{AttCap} and our
assembly map is the blocked surgery leafwise assembly map. The cosheaf 
$\mathcal{L}^{uf}_*(S_x)$ assigns to each simplex $\sigma$ of 
$BG$, the group $H^{uf}_*(S_\sigma;\textbf{L})$, which is the uniformly 
finite L-homology of the union of the leaves of 
$\mathcal{F}$ corresponding to $\sigma$ in $BG$. 
Taking the sum of the simplices of $BG$ with coefficients in 
$\mathcal{L}^{uf}(S_x)$ as in \cite{Ran} chapter 4,  we get
the normal invariants on the manifold $M$. We let $\mathcal{A}$ be the additive
category of free modules over the cosheaf $\mathcal{L}^{uf}$ and 
the simplicial complex $K$ is the classifying space $BG$.

In this case $\mathcal{A}_*[K]$ is the category of functors that assigns to
a simplex $\sigma \in BG$ the group $\mathcal{L}_*^{uf}(S_\sigma)$, the L-homology of
the union of the leaves corresponding to $\sigma$. We then have the category
$\mathcal{A}_*(K)$ as the category of direct sums 
$$\sum_{\sigma \in BG} \mathcal{L}_*^{uf}(S_\sigma)$$
We then have the covariant assembly functor \cite{Ran}
$$\mathcal{B}(\mathcal{A}_*[K]) \to \mathcal{B}((\mathcal{A})^*(K))$$
the reason for the upper * being that the chain complex $\Delta(K)$ of $K$
is taken to $\Delta_*[K](\sigma)=S^{\mid\sigma\mid}\Delta(\sigma)$ and 
therefore ends up in homology.
The assembly map  takes a chain complex $C_*(BG;\mathcal{L}_{uf})=
C^{uf}_*(BG;\mathcal{L}_{uf})$ (since $BG$ is compact) to the 
complex
$$C_*[BG]_r=\sum_{\sigma \in BG}C_{r-\mid \sigma \mid}(\sigma; \mathcal{L}^{uf}(S_{\sigma})),$$
where the dimension shift $r-\mid \sigma \mid$ is 
the chain duality $T$ of Definition 3.9 (see Example 4.13 of \cite{Ran} and 
Definition 3.9 above).
We then use the Leray-Serre spectral sequence of a cosheaf \cite{Ran}, 
\cite{Quinn}
$$H_*(\sum_{\sigma \in BG}C_{r-\mid \sigma \mid}(\sigma; \mathcal{L}^{uf}(S_\sigma)))\Rightarrow H_*(M;\textbf{L})$$
which converges to $H_*^{uf}(M;\textbf{L})=H_*(M;\textbf{L})$ as $M$ is 
compact and clearly yields the assembly map to $H_*(M;\textbf{L})$. 

This gives us an assembly map 
$$H_*(BG;\mathcal{L}^{uf}(S_x)) \to H_*(M;\textbf{L})$$
to the L-homology of $M$.
Let $\Delta^{(\infty)}(\nu(S_x)) \in KO^{uf}_*(S_x) \otimes \textbf{Q}$, be 
the uniformly finite $KO[1/2]$-orientation of the normal bundle of $\nu(S_x)$
(See \cite{Ran}, Chapter 16). 
We further have maps 
$$ph:KO^*(M)\otimes \textbf{Q} \to H^*(M;\textbf{Q})$$ 
given by the ordinary Pontrjagin character, $ph(\xi)=ch(\xi \otimes \textbf{C})$,
where $ch$ is the Chern character. We also have the $L^\infty$ Pontrjagin 
character
$$ph^{(\infty)}:KO_{uf}^*(S_x)\otimes \textbf{Q} \to H^*_{uff}(S_x; \textbf{Q})$$
given by the $ph^{(\infty)}(\xi)=ch^{(\infty)}(\xi \otimes \textbf{C})$,
where $ch^{(\infty)}$ is the $L^\infty$ Chern character defined as follows.
Let $(E,\nabla)$ be a vector bundle of bounded geometry over $S_x$. Suppose
it is split into line bundles
$$(E,\nabla)=(L_1,\nabla_1)\oplus ... \oplus (L_n, \nabla_n)$$
Then we have 
$$ch^{(\infty)}(E,\nabla)=e^{c_1^{(\infty)}(L_1,\nabla_1)}+...+e^{c_1^{(\infty)}(L_n, \nabla_n)}$$ 
where $c_1^{(\infty)}$ is the $L^\infty$ first Chern class.
In analogy with the Atiyah-Singer index theorem for families, we have a 
map $x \mapsto [\Delta^{(\infty)}(\nu(S_x))]$ assigning to each leaf in $BG$ its 
$KO[1/2]$-orientation as a K-theory class in $KO^{uf}_*(S_x)\otimes \textbf{Q}$. 
The assembly map takes the 
homology class $\sum [\Delta^{(\infty)}(\nu(S_x))]x$ to the $KO[1/2]$-orientation 
$[\Delta(\nu(\mathcal{F}))]$ in $KO(M)\otimes \textbf{Q}$. 
The image of this under the  assembly map is therefore
 $[\Delta(\nu(\mathcal{F}))]=\Delta([TM/E])$.
Computing the images on both sides of $ph$ and $ph^{(\infty)}$, we 
obtain the result that the Hirzebruch L-classes $L_0^{(\infty)},
L_1^{(\infty)},...,L_k^{(\infty)}$ of $\nu(S_x)$ assemble to the L-classes
$L_0,L_1,...,L_k$ of $\nu(\mathcal{F})$. 
We then have that the forms $L_0^{(\infty)},L^{(\infty)}_1,...,L^{(\infty)}_k$
generate the Pontrjagin ring $Pont^{(\infty),q}(\nu(S_x))$ and that
$L^{(\infty)}_q=0$ if $q>k$ because of the dimension of the normal bundle of
the leaf which is $k$ where $k$ is the codimension of the foliation.
From this, we have that the forms $L_i$ satisfy $L_q=0$ if $q>k$ since
the $L^{(\infty)}_i$ assemble to them.
We then have the classes $p_0, p_1,...,p_k$, with $p_q=0$ if $q>k$,
 which generate $Pont^*(TM/E)$ and the
result follows by Shulman's theorem below\cite{BSS}.  
\begin{thm}[Shulman \cite{BSS}] Let $G$ be a Lie group.
A real characteristic class $\Phi$ associated to an invariant polynomial of 
degree $q$ on the Lie algebra of $G$, has a representative in the de Rham 
complex that involves only forms of degree $\ge q$. 
Thus:
$$\Phi=\Phi_{q-1}+\Phi_{q-2}+...+\Phi_0$$
where $\Phi_i \in \Omega^{q+i}N_{q-i}G$, and $\Omega$ is the de Rham functor
and $\Omega N G$ is the de Rham complex of $NG$. Here $NG$ is the 
Eilenberg-MacLane bar construction
$$NG: pt \gets G \mbox{ }\substack{\gets\\[-1em] \gets}\mbox{ } G \times G \mbox{ }\substack{\gets\\[-1em]\gets\\[-1em] \gets}\mbox{ } G \times G \times G \mbox{ }\substack{\gets\\[-1em]\gets\\[-1em]\gets\\[-1em]\gets} \mbox{ }... $$
\end{thm}
Here $\Omega^q(M)$ denotes $\Gamma\bigwedge^qT^*(M)$, where $\bigwedge^iT^*$
denotes the i-th exterior power of the cotangent bundle of $M$ and $\Gamma$
is the $C^\infty$-sections of a vector bundle \cite{Bott1}.

The double complex $\Omega NG$
is has two differentials: a de Rham differential $d:\Omega^i NG \to
\Omega^{i+1}NG$ and a differential operator
$$\delta:\Omega^q N_p(G)\to\Omega^q N_{p+1}(G)$$
which commutes with $d$. There is a total operator $D=d\pm\delta$ \cite{BHoch}.

If one filters by $\le p$, the resulting spectral sequence obviously has the
$E_1$-term:
$$E_1^{p,q}=H^q(N_pG)$$
and $d_1$, induced by $\delta$ yields the cobar construction on $H^*(G)$.
It follows that the $E_2$ term of this sequence is the symmetric algebra
generated by the primitive elements $\mathcal{P}(G)$ of $H^*(G)$, considered
in $H^{1,*}(G)$ and all further differential operators vanish.
This leads to the evaluation of $HNG$ \cite{BHoch}:
$$HN(G)\simeq S\mathcal{P}(G)$$
\par\noindent
\textit{Sketch of Proof:} \cite{BHoch} Let $\mathcal{G}$ be the Lie algebra of $G$. The Chern-Weil
construction defines a homomorphism from the invariants $Inv_G(S\mathcal{G}^*)$
of $S\mathcal{G}^*$ to $H^*(BG)$. We have a map
$$\phi:Inv_G(S\mathcal{G}^*) \to H^*(BG)$$
The cohomology of the double complex is:
$$H_\delta^p(\Omega^q NG) \simeq H_{cont}^{p-q}(G; S^q \mathcal{G}^*)$$
This is proven in \cite{BHoch}. We have
$$H_{cont}^0(G;W)=Inv_G W$$
and 
$$H_{cont}^i(G;W)=0 \mbox{ for } i>0$$
if $G$ is compact. Thus
$$H_\delta^p\Omega^p G=Inv_G S^p \mathcal{G}^*$$
whence, as there is no cohomology above the diagonal, the ``edge-homomorphism''
induces a natural map
$$\Phi:Inv_G(S\mathcal{G}^*) \simeq H^*(BG)$$
which is the Chern-Weil homomorphism.
Thus we have a cocycle $\Phi_0$ in $\Omega NG$, which can be completed by
adding correction terms $\Phi_i$, $i>0$ because the $\delta$-cohomology
vanishes.

This proves Bott 
integrability in the case where there is no fundamental group.
\begin{rk}
The proof above works exactly for a topological microbundle $E$ and a 
topological foliation as for a smooth vector bundle and a smooth foliation,
except for Shulman's theorem which depends on the Pontrjagin classes 
being expressed in terms of a curvature form. Tsuboi \cite{Tsuboi, Hurder}
has shown the Haefliger space to be homotopy equivalent to $BO_q$ in the
case of $C^1$ foliations, showing that a topological version of Bott
Integrability is false.
\end{rk}
We also have the corresponding theorem for Chern classes.
\begin{dfn} Let $E$ be a vector bundle over a complex manifold $M$. 
The Chern ring of E, $Chern^*(E)$ is the ring in $H^*(M;\textbf{C})$ generated
by the Chern classes of E.
\end{dfn}
A parallel result for Chern classes of complex manifolds was obtained by Bott.
\begin{thm} \cite{Bott} Let $M$ be a complex manifold, $E$ a holomorphic subbundle of 
$TM$ which is integrable. Let $Chern^*(TM/E)$ be the Chern ring of $TM/E$. Then
$Chern^q(TM/E)=0$ for $q > 2k$. 
\end{thm}
The proof is as before except with $H^{uf}_n(M;\textbf{L}(\textbf{C}))$
in place $H^{uf}_n(M;\textbf{L})$, and $K(M)$ and $\mathcal{K}_{uf}(S_x)$ in 
place of $KO(M)$ and $\mathcal{KO}_{uf}(S_x)$, along with the Chern 
character in place of the Pontrjagin character.

There is an alternative approach by Crainic and Moerdijk \cite{CrMo}.
\begin{dfn}\cite{CrMo} Let $M$ be a manifold of dimension $n$ equipped with a foliation
$\mathcal{F}$ of codimension $q$. A transversal section of $\mathcal{F}$ is an
embedded $q$-dimensional submanifold $U \subset M$ which is everywhere 
transverse to the leaves. Recall that if $\alpha$ is a path between
two $x$ and $y$ on the same leaf, and if $U$ and $V$ are transversal 
sections through $x$ and $y$, then $\alpha$ defines a transport along the 
leaves from a neighborhood of $x$ in $U$ to a neighborhood of $y$ in $V$. 
If transport ``along $\alpha$'' is defined in all of $U$ and embeds $U$ 
into $V$, it will be called a holonomy embedding. 
\end{dfn}
\begin{dfn}\cite{CrMo}
Transversal sections $U$ through $x$ as above should be thought of as 
neighborhoods of the leaf through $x$ in the leaf space. This motivates
the definition of a transversal basis for $(M,\mathcal{F})$ as a family
$\mathcal{U}$ of transversal sections $U \subset M$ with the property
that, if $V$ is any transversal section through a given point $y \in M$
there exists a holonomy embedding $h:U \to V$ with $U\in\mathcal{U}$ and
$y\in h(U)$.
\end{dfn}
\begin{dfn}\cite{CrMo}
Typically, a transversal section is a $q$-disk given by a chart for the
foliation. Accordingly, we can construct a transversal basis $\mathcal{U}$
out of a basis $\tilde{\mathcal{U}}$ of $M$ by domains of foliation charts
$\phi_U:\tilde{U}\to \textbf{R}^{n-q} \times U$, 
$\tilde{U}\in\tilde{\mathcal{U}}$ with $U=\textbf{R}^q$. Note that each 
inclusion $\tilde{U}\to\tilde{V}$ between opens of $\tilde{\mathcal{U}}$
induces a holonomy embedding $h_{U,V}:U \to V$ defined by the condition
that the plaque in $\tilde{U}$ through $x$ contained in the plaque in 
$\tilde{V}$ through $h_{U,V}(x)$. 
\end{dfn}
\begin{dfn}\cite{CrMo} Let $(M, \mathcal{F})$ be a foliated manifold and let 
$\mathcal{U}$ be a transversal basis. Consider the double complex
which in bidegree $k,l$ is the vector space
$$C^{k,l}=\check{C}^k(\mathcal{U},\Omega^l)=\prod_{U_0\overset{h_1}{\to}...
\overset{h_k}{\to}U_k}\Omega^l(U_0)$$
Here the product ranges over all $k$-tuples of holonomy embeddings between
transversal sections from the given basis $\mathcal{U}$, and $\Omega^k(U_0)$
is the space of differential $k$-forms on $U_0$. For elements $\omega \in 
C^{k,l}$, we denote its components by 
$\omega(h_1,...,h_k)\in\Omega^k(\mathcal{U}_0)$. The vertical differential
$C^{k,l}\to C^{k,l+1}$ is $(-1)^kd$ where $d$ is the usual De Rham 
differential. The horizontal differential $C^{k,l}\to C^{k+1,l}$ is 
$\delta=\sum (-1)^i\delta_i$ where
$$\delta_i(h_1,...,h_{k+1})=\begin{cases} h_1^*\omega(h_2,...,h_{k+1})
\mbox{ if }i=0\\\omega(h_1,...,h_{i+1}h_i,...,h_{k+1})\mbox{ if }
0<i<k+1\\\omega(h_1,...,h_k)\mbox{ if }i=k+1\end{cases}$$
This double complex is actually a bigraded differential algebra, with the
usual product
$$(\omega \cdot \eta)(h_1,...,h_{k+k^\prime})=(-1)^{kk^\prime}\omega(h_1,...
h_k)h_1^*...h_k^*\eta(h_{k+1},...,h_{k+k^\prime})$$
for $\omega \in C^{k,l}$ and $\eta \in C^{k^\prime, l^\prime}$. We will also
write $\check{C}(\mathcal{U},\Omega)$ for the associated total complex,
and refer to it as the $\check{C}$ech-DeRham complex of the foliation. 
The associated cohomology is denoted $\check{H}^*_{\mathcal{U}}(M/\mathcal{F})$
and referred to as the $\check{C}$ech-DeRham cohomology of the leaf space 
$M/\mathcal{F}$ with respect to the cover $\mathcal{U}$.
\end{dfn}
\begin{dfn}\cite{CrMo}
In general, choosing a transversal basis $\mathcal{U}$ and a basis
$\tilde{\mathcal{U}}$ of $M$ there is an obvious map of double complexes
$C^{k,l}(\mathcal{U})\to C^{k,l}(\tilde{\mathcal{U}})$ into the
$\check{C}$ech-DeRham complex of the manifold $M$. Hence a canonical 
map
$$\pi^*:\check{H}_{\mathcal{U}}^*(M/\mathcal{F})\to H^*(M;\textbf{R})$$
which should be thought of as the pullback along the ``quotient map''
$\pi:M \to M/\mathcal{F}$.
\end{dfn}
\begin{thm}[Crainic-Moerdijk]\cite{CrMo} There is a natural isomorphism
$$\check{H}_{\mathcal{U}}^*(M/\mathcal{F})\simeq H^*(BG;\textbf{R})$$
between the $\check{C}$ech-DeRham cohomology and the cohomology of the
classifying space of the holonomy groupoid.
\end{thm}
\par\noindent\textit{Alternate Proof of Bott Integrability:}
Now let us note that $H_*(BG;\textbf{L})\otimes \textbf{R}
=H_*(BG;\textbf{R})$, where $\textbf{L}=G/TOP$ is the simply-connected 
surgery spectrum. Similarly $H_*(M;\textbf{L})\otimes \textbf{R}=
H_*(M;\textbf{R})$. We can then use the map 
$$\pi^*:H^*(BG;\textbf{R})\to H^*(M;\textbf{R})$$
 to pull back the $L$-classes from $BG$ to $M$. In 
addition to this the $L$-classes of the tangent bundle of $BG$ vanish above 
dimension $q$ if $q$ is the codimension of the foliation from the corresponding
 fact for the normal bundles to the leaves $L_\sigma$ corresponding to the 
simplex $\sigma \in BG$ and from the spectral sequence of a cosheaf above:
$$H_*(\sum_{\sigma \in BG}C_{r-\mid \sigma \mid}(\sigma; \mathcal{L}^{uf}(S_\sigma)))\Rightarrow H_*(M;\textbf{L})$$
 we get a map to $H_*(M;\textbf{L})$. 
By Shulman's theorem, the Pontrjagin ring vanishes above dimension 
$2q$. Hence if we pull back the Pontrjagin ring to $H^*(M;\textbf{R})$ we
get that the Pontrjagin ring of $TM/E$ vanishes above dimension $2q$.
\section{Higher Integrability}
\begin{prop} The Novikov conjecture holds for a group $\Gamma$ if and only
  if the assembly map in symmetric $L$-theory
  $$A_\Gamma:H_*(B\Gamma;\textbf{L}(\textbf{Z}))\to L^*(\textbf{Z}[\Gamma])$$
  is a rational split injection. Note that $L^*(\textbf{Z}[\Gamma])\otimes
  \textbf{Z}[1/2]=L_*(\textbf{Z}[\Gamma])\otimes\textbf{Z}[1/2]$, where
  $L^*$ is the symmetric $L$-theory and $L_*$ is ordinary $L$-theory
  \cite{FRR}.
\end{prop}
\begin{dfn} The Atiyah-Hirzebruch spectral sequence for uniformly finite
  homology is the spectral sequence
  $$E_{p,q}^2=H_p^{uf}(X;\pi_q(\textbf{E}))\Rightarrow E_{p,q}^\infty=H_{p+q}^{uf}(X;\textbf{E})$$
  whose $E_2$ term is the uniformly finite homology with coefficients in the
  homotopy groups of the spectrum $\textbf{E}$ with norm given by the metric
  on $\textbf{E}$ and $E_\infty$ term given by the generalized uniformly
  finite homology with coefficients in the spectrum $\textbf{E}$.
  \end{dfn}
\begin{dfn} Suppose now we have a fundamental group $\pi_1(M)$ for $M$,
and that $M$ admits a foliation $\mathcal{F}$ with classifying space for the
holonomy groupoid $BG$. Let $\pi$ be a group, and let $j:\pi_1(BG) \to \pi$ be
 a homomorphism. We can consider the corresponding 
 map $i:M \to B\pi$ and let $y \in H^i(B\pi; \textbf{Z})$, with $i^*(y)=\gamma$
 then we have higher Pontrjagin classes $\gamma \cup p_i(TM/E)$ where
 $E \subset TM$ is a subbundle of the tangent bundle of $M$. Consider the
 bundle $E$ as above, and suppose it is integrable with foliation
 $\mathcal{F}$. We have the higher Pontrjagin ring $Pont^*(\pi)(TM/E)$
 consisting of all products with higher
Pontrjagin classes $\gamma \cup p_i(E)$ for $\gamma=i^*(y)$, $y \in H^*(B\pi; \textbf{Z})$.
\end{dfn}
\begin{dfn} We have a dual formulation in terms of $L$-classes in 
homology. Let $M$ be a compact manifold of dimension $n$. 
Suppose $E \subset TM$ is integrable, $\mathcal{F}$ the corresponding
foliation and let $BG$ be the classifying space of the foliation groupoid,
and $j:\pi_1(BG)\to \pi$ be a homomorphism. Then the higher Pontrjagin ring in $H_*(M)$ is defined
as follows. Let $I^*:H^*(B\pi) \to H^*(M)$ be the homomorphism induced from 
the map $M \to B\pi$. Let $L_i(TM/E)$ be the $i$-th $L$-class in the homology
of $M$ of $TM/E$, the Poincar\'e dual of the $L$-classes in cohmology. We 
then have classes of the form $I^*(\gamma)\cap L_i(TM/E)$ which generate 
a ring in homology under the intersection product, which is the 
Poincar\'e dual of the cup product.
\end{dfn}
\begin{lmm}[Higher Shulman's Theorem] Let $(M,\mathcal{F})$ be a foliated
  manifold and let $BG$ be the leafspace of $\mathcal{F}$. Let $\pi$ be
  the fundamental group of $BG$. If the higher Pontrjagin classes of
  $(M, \mathcal{F})$ all vanish above the codimension $k$ of $\mathcal{F}$,
  then the higher Pontrjagin ring $Pont(\pi)^q(\nu(\mathcal{F}))$ vanishes
  above dimension $2k$, where $\nu(\mathcal{F})=TM/E$, $E$ being the
  tangent bundle of $\mathcal{F}$.
\end{lmm}
\textit{Proof:} Note that we have a map $BG \to B\pi$. We consider the
cohomology of $BG$ with coefficients in the cosheaf of uniformly finite
$\textbf{L}$-cohomology groups of the transversals, where $\textbf{L}$ is
the surgery spectrum:
$$H^*(BG; \mathcal{L}_{uf}^*(T_x))$$
where $T_x$ is the transversal through $S_x$ at $x \in BG$ the leaf through
$x \in BG$ and $\mathcal{L}^*_{uf}$ is the $\textbf{L}$-cohomology
cosheaf. From the standard Leray-Serre spectral sequence we get an
assembly map
$$H^*(BG;\mathcal{L}^*_{uf}(T_x)) \to H^*(M;\textbf{L})$$
We then tensor with the reals:
$$H^*(BG;\mathcal{L}^*_{uf}(T_x)) \otimes \textbf{R} \to H^*(M;\textbf{L})\otimes \textbf{R}$$
We then use the fact that $\textbf{L}\otimes\textbf{Q}=\textbf{KO}\otimes\textbf{Q}$
and we have
$$H^*(BG; \mathcal{KO}^*_{uf}(T_x))\otimes\textbf{R}\to H^*(M;\textbf{KO})\otimes\textbf{R}$$
where $\mathcal{KO}^*_{uf}(T_x)$ is the cosheaf of uniformly finite
KO-cohomology groups of the transversals. We then apply the Pontrjagin
character to get
$$H^*(BG; \mathcal{H}_\beta^*(T_x))\otimes\textbf{R}\to H^*(M;\textbf{R})$$
where $\mathcal{H}_\beta^*(T_x)$ is the cosheaf of bounded de Rham cohomology
of the transversals. We observe that the de Rham forms representing elements
of the $L^\infty$ Pontrjagin ring of $T_x$ have a decomposition via
Shulman's theorem:
$$\Phi=\Phi_{q-1}+\Phi_{q-2}+...+\Phi_0$$
where $\Phi_i \in \Omega^{q+1}N_{q-i}SO(n)$ and $\Omega$ is the de Rham
functor and $\Omega NSO(n)$ is the de Rham complex of $NSO(n)$. We then
use the map $BG \to B\pi$ to get the assembly map
$$H^*(B\pi; \mathcal{H}_\beta^*(T_x))\otimes\textbf{R}\to H^*(M;\textbf{R})$$
Elements of the higher Pontrjagin ring of $T_x$ can be represented by
elements to be represented as a de Rham form $\gamma \otimes \Phi$, where
$\gamma$ is a de Rham form representing an element of $H^*(B\pi)$ and
$\Phi$ is a representative of the Pontrjagin ring of $T_x$. But
after assembling to $M$ we see that this yields a de Rham form representing
an element of the higher Pontrjagin ring of $TM/E$ and that this is zero
above dimension $2k$.

\begin{thm} Let $M$ be a compact closed manifold of dimension $n$, suppose 
$E \subset TM$ is a subbundle of the tangent bundle of $M$. Let $E$ be 
integrable, $BG$ the classifying space of the foliation groupoid of the 
corresponding foliation $\mathcal{F}$ and $j:\pi_1(BG)\to \pi$ a homomorphism.
Suppose the Novikov conjecture 
is true for $\pi$. Then the higher Pontrjagin ring
$Pont^q(\pi)(TM/E)$ vanishes above dimension $2k$, where $k$ is the
dimension of $TM/E$.
\end{thm}
\textit{Proof:} Let $k=dim(TM/E)$, $n=dim(TM)$. We start with an assembly map 
as before. This assembly map is analogous to the one in blocked surgery
\cite{Quinn, CW}. Let $BG$ be the classifying space of the foliation groupoid. 
We take the homology of $BG$ with coefficients in the cosheaf 
$\mathcal{L}^{uf}_*(S_x)$ (for homology with coefficients in a cosheaf of 
spectra, see the stratified surgery exact sequence in \cite{W}), 
where $\mathcal{L}^{uf}(S_x)$ is the 
uniformly finite $\textbf{L}(\pi)$-homology of $S_x$, the leaf through $x \in BG$, where the classifying space of the holonomy groupoid $BG$ is considered
as the leaf space of $\mathcal{F}$ (to identify $BG$ with the leaf space of $\mathcal{F}$, see \cite{Moerd}),
and $\textbf{L}(\pi)$ is the surgery spectrum of $\pi$.
We get an assembly map
$$H_*(BG ; \mathcal{L}^{uf}(S_x)) \to H_*(M;\textbf{L}(\pi))$$
where $\textbf{L}(\pi)$ is the surgery spectrum of $\pi$.
We next define the higher signatures of the leaves. Let $i:M \to B\pi$, and
$j:S_x \to M$ be the embedding of the leaf $S_x$ through $x \in BG$ so that
for $y \in H^*(B\pi), \gamma=j^*i^*(y)\in H_{uf}^*(S_x)$. Let
$\Delta^{(\infty)}(\nu(S_x))\in KO_*^{uf}(S_x)\otimes \textbf{Q}$ be the
$KO[1/2]$ orientation of the normal bundle $\nu(S_x)$ of the leaf $S_x$
through $x \in BG$. Define $\beta=ph^{(\infty)-1}(\gamma)$. We can then
define the higher signature class $\Delta_\gamma^{(\infty)}(\nu(S_x))=
\beta \cup \Delta^{(\infty)}(\nu(S_x)) \in H_*^{uf}(S_x;\textbf{KO}(\pi))
\otimes \textbf{Q}$ be the higher signature of the normal bundle
$\nu(S_x)$ to the leaf $S_x$ through $x, x\in BG$. This higher signature
class is detected by the surgery space $\textbf{L}(\pi)$. This class lies
in $H_*^{uf}(S_x; H_*(B\pi; \textbf{L}))$ and by the Novikov conjecture
this group injects into $H_*^{uf}(S_x; \textbf{L}_*(\pi))$. This is because
the symmetric assembly map is a split injection on ordinary homology, and
hence we can use the Zeeman Comparison Theorem \cite{Z}, which states that
if the $E_2$ terms of two spectral sequences are isomorphic, then so are
the $E_\infty$ terms, on the summand to obtain an injection.

We have the following commutative diagram:
$$\begin{CD}
  H_*(BG;\mathcal{H}^{uf}(\pi)_*(S_x)\otimes\textbf{Q})@<ph^{(\infty)}<< H_*(BG;\mathcal{KO}^{uf}(\pi)(S_x))\otimes\textbf{Q}@>A^{(\infty)}_{Surgery}>> H_*(BG;\mathcal{L}^{uf}(S_x))\otimes\textbf{Q}\\
  @VA_{Leaves}^{\mathcal{H}}VV @VA_{Leaves}^{\mathcal{KO}}VV @VA_{Leaves}^{\mathcal{L}}VV \\
  H_*(M; H_*(B\pi; \textbf{Q})) @<ph<< H_*(M;\textbf{KO}_*(\pi)) \otimes \textbf{Q} @>A_{Surgery}>> H_*(M;\textbf{L}(\pi))\otimes \textbf{Q}
\end{CD}$$
Here $ph$ is an isomorphism, and $A_{Surgery}:H_*(M;\textbf{KO}_*(\pi))\otimes\textbf{Q} \to
H_*(M;\textbf{L}(\pi))\otimes\textbf{Q}$ is injective as $KO_*(\pi)\otimes \textbf{Q}\to L(\pi)\otimes \textbf{Q}$
is by the Novikov conjecture for $\pi$. 

In addition we have the cosheaf $\mathcal{L}^{uf}(S_x)$
as above, and the other cosheaves $\mathcal{H}^{uf}(\pi)(S_x)$ which associates
to $x$ the uniformly finite $B\pi$-homology of the leaf through $x$ and
$\mathcal{KO}^{uf}(\pi)_*(S_x)$ is the cosheaf which associates to $x$ the
uniformly finite $KO(\pi)$-homology of the leaf through $x$.
The map
$$A_{Surgery}^{(\infty)}:H_*(BG;\mathcal{KO}^{uf}(\pi)(S_x))\otimes\textbf{Q}\to H_*(BG;\mathcal{L}^{uf}(S_x))\otimes\textbf{Q}$$
is the assembly map coming from the injective map $KO(\pi)\otimes\textbf{Q}\to L(\pi)\otimes\textbf{Q}$ and is
itself injective.
The maps $A_{Leaves}^{\mathcal{H}}, A_{Leaves}^{\mathcal{KO}},
A_{Leaves}^{\mathcal{L}}$ are assembly maps of the cosheaves.

We have $\sum_x \Delta^{(\infty)}_\gamma(\nu(S_x))[x]$ as a cycle in generalized homology
\cite{Ran}, Chap. 12.
These each assemble to a higher signature class 
$\Delta_\gamma(\nu(\mathcal{F}))$, with $\gamma \in H^*(B\pi)$. 
We then have classes in $H_*(M;L_*(\pi))$ coming from the 
leaves. We take the $L^\infty$ Pontrjagin 
character of $\Delta^{(\infty)}_\gamma(\nu(S_x))$, and obtain the higher
total L-class $\gamma \cup \mathcal{L}_{Tot}^{(\infty)}(\nu(S_x))=
\gamma\cup(1+L_1^{(\infty)}(\nu(S_x))+L_2^{(\infty)}(\nu(S_x))+...)$,
which yields $\gamma \cup L^{(\infty)}_0
(\nu(S_x))=\gamma,\gamma\cup L^{(\infty)}_1(\nu(S_x)),...,
\gamma\cup L^{(\infty)}_k(\nu(S_x))$. Note that the Pontrjagin character
is a rational isomorphism. For dimensional reasons, these higher
$L$-classes vanish above the dimension of $\nu(S_x)$ which is $k$.

Because of the Novikov conjecture, if the class of the surgery obstruction
$\Delta_\gamma^{(\infty),1}(\nu(S_x)) \in H_*^{uf}(S_x;\textbf{L}(\pi))$
coming from the higher signature \break
$ph^{(\infty)}(\Delta_\gamma^{(\infty)}(\nu(S_x)) \in H_*^{uf}(S_x; KO(\pi))$ is zero then the higher signature is zero. So the vanishing of the
higher Pontrjagin classes from $KO(\pi)$ for dimensional reasons guarantees by the
Novikov conjecture that the higher Pontrjagin classes from the surgery group
vanish in the same range.
These classes assemble via $A_{Leaves}^{\mathcal{L}}$ to $\Delta^1_\gamma(\nu(\mathcal{F}))$ in
$H_*(M;\textbf{L}(\pi))$. Again we have a higher signature class
$\Delta_\gamma(\nu(\mathcal{F}))$ in $H_*(M;KO(\pi))$, which injects
into $H_*(M;\textbf{L}(\pi))$ via $A_{Surgery}$ because of the Novikov
conjecture for $\pi$.

The Pontrjagin character of $\Delta_\gamma(\nu(\mathcal{F}))$ is of the
form $ph(\beta \cap \Delta(\nu(\mathcal{F})))$ where $\Delta(\nu(\mathcal{F}))$
is the KO[1/2]-orientation of $\nu(\mathcal{F})$, and $\beta \in KO(M)$ is
such $ph(\beta)=\gamma$, representing $\Delta_\gamma(\nu(\mathcal{F}))$
in $H_*(M;KO(\pi))$ and hence comes from an invariant polynomial in the
curvature form of $\Delta_\gamma(\nu(\mathcal{F}))$. This comes from the definition of
the Chern character:
$$ch(\xi)=\sum_{j \in \textbf{N}}k_j tr(F_\nabla \wedge ... \wedge F_\nabla)$$
as a de Rham form, with $\xi$ a vector bundle, 
where $\nabla$ is a connection on $\xi$, $k_j=\frac{1}{j!}(\frac{1}{2\pi i})^j$
and $F_\nabla$ is the curvature of this connection.
The Pontrjagin character is then defined as
$ph(\xi)=ch(\xi \otimes \textbf{C})$.

Because of the vanishing over dimension $k$ of the $L^\infty$ higher
Pontrjagin classes of the leaves we have vanishing of the higher Pontrjagin
classes of
$\nu(\mathcal{F})$ above dimension $k$, from the assembly map.
We have vanishing over dimension $k$ of the higher Pontrjagin classes and
hence vanishing over dimension $2k$ of the higher Pontrjagin ring
after applying Lemma 4.1 (Higher Shulman's theorem).
\par\noindent\textit{Alternate Proof:}
We have an alternative approach to this theorem using Theorem 3.5 above, the
definitions which precede it and \cite{CrMo}. Let us 
consider $H^*(BG;\textbf{L}(\pi))$ where there is a homomorphism
$j:\pi_1(BG) \to \pi$. Using the 
Novikov conjecture for $\pi$ we get an injection $KO(\pi)\otimes\textbf{Q}\to \textbf{L}(\pi)\otimes\textbf{Q}$.
We thus have an injection
$H^*(BG;KO(\pi)\otimes \textbf{R})\to H^*(BG;\textbf{L}(\pi)) \otimes\textbf{R}$.
As before this is a split injection by Proposition
4.1, and so the higher signatures, and therefore the L-classes are detected by
$\textbf{L}(\pi)$. Now apply the Pontrjagin character
$$ph \otimes \textbf{R}: H^*(BG;KO(\pi) \otimes \textbf{R}) \to
H^*(BG;H^*(B\pi;\textbf{R}))$$
So if we consider the map
$\rho^*:H^*(BG;H^*(B\pi;\textbf{R}))\to H^*(M;H^*(B\pi;\textbf{R}))$
from Theorem 3.5, 
we can pull back the higher Pontrjagin classes of $BG$ to $M$.
Since the higher Pontrjagin classes vanish above the
dimension $k$ of the tangent bundle of the leafspace $BG$, $k$ being
the codimension of the foliation, from the corresponding fact for
the normal bundles of the leaves, the Pontrjagin classes of $TM/E$ vanish
above $k$ from being in the image of the map $\rho^*$.
We now apply Lemma 4.1 (Higher Shulman's Theorem) to finish the proof. 

\begin{cor} Suppose $M$ is a $K(\pi,1)$ with $\pi$ satisfying the Novikov
  conjecture and that $E$ is a subbundle of the
  tangent bundle with $p_i(TM/E)\ne 0$ for some $i \ge 0$ and $p_i$ a real
  Pontrjagin class. Suppose in addition
  that $p_i(TM/E) \cup c \ne 0$ for some class $c \in H^*(M; \textbf{R})$
  and dim$(p_i \cup c)>$2dim(TM/E).
  Then $E$ cannot be the tangent bundle to a foliation with simply connected
  leaves.
\end{cor}
\begin{rk} The reason for the hypothesis that the leaves of the foliation
  be simply connected is that $\pi_1(M)=\pi_1(BG)$ in this case. This follows
  from Corollary 3.2.4 of Haefliger \cite{Haefliger}:
  \par\noindent
Let $\mathcal{F}$ be a foliation on a manifold $X$ so that the holonomy
coverings of the leaves are $(k-1)$-connected. Then the holonomy groupoid
$G$ of $\mathcal{F}$ considered as a $G$-principal bundle with base $X$ by
the end projection is $k$-universal. Thus the space $X$ itself is
$k$-classifying and the map $i$ of $X$ to $BG$ is $k$-connected.

Corollary 3.1.5 of Haefliger \cite{Haefliger} states the following:
\par\noindent
Suppose the target projection of the holonomy groupoid $G$ of the foliation
$\mathcal{F}$ on a manifold $X$ is a locally trivial fibration, whose fiber is
the common holonomy covering $L$ of all the leaves. Then the map $i:X \to BG$
is homotopy equivalent to a locally trivial fibration with base $BG$
and fiber $L$.

The hypothesis of Corollary 3.1.5 of Haefliger are satisfied in the following
cases \cite{Haefliger}:
\par\noindent i. $X$ is compact and $\mathcal{F}$ possesses a transverse
Riemannian metric.
\par\noindent ii. The leaves of $\mathcal{F}$ are transverse to the fibers
of a compact fibration and the foliation is analytic.
\par\noindent iii. The leaves of $\mathcal{F}$ are the trajectories of a flow
without closed orbits or the holonomy group of each closed orbit is infinite.

We can use these results to do the following. Consider a leaf $i:L \to X$.
We can take the covering corresponding to $i_*(\pi_1(L))$ of $\tilde{X}$ of
$X$ and the induced foliation $\tilde{\mathcal{F}}$ of $\mathcal{F}$. If
$L$ is the common holonomy covering of all the leaves, then Haefliger's
result means that $\pi_1(X)=\pi_1(BG)/i_*(\pi_1(L))$. 
\end{rk}
\begin{cor} Suppose the leaves of $\mathcal{F}$ are all simply connected.
  Then $\pi_1(M)=\pi_1(BG)=\pi$. Suppose the Novikov conjecture is true
  for $\pi$. Then the higher Pontrjagin ring of $TM/E$, $Pont^q(\pi)(TM/E)$
  vanishes for $q>2k$ where $k=dim(TM/E)$.
  \end{cor}
\begin{cor} Suppose $M$ is the hyperbolization \cite{DJ} of a manifold $N$.
  Then the results of the previous corollary hold.
\end{cor}
\begin{cor} Let $\mathcal{F}$ be a foliation with contractible hyperbolic leaves. Then
  $Pont^{k}(\pi)(TM/E)=0$ if $k > 2$, provided the Novikov Conjecture is
  true for $\pi=\pi_1(M)$.
\end{cor}
\textit{Proof:} The uniformly finite cohomology of any leaf
$H_{uf}^k(L;\textbf{Q})=0$ for $k>1$ \cite{Att} and \cite{Gromov}, Example
0.1.B.
So following through the proof of Theorem 4.1 we have that
$Pont^{k}(\pi)(TM/E)=0$ if $k>2$, by Shulman's theorem
with $q=1$. Since the homotopy
type of $M$ is that of $BG$, by Haefliger's Corollary 3.2.4,
the result follows.
\begin{cor} Let $\mathcal{F}$ be a foliation with all leaves a contractible
  symmetric space of rank $m$. Then $Pont^{k}(\pi)(TM/E)=0$ if $k>2m$, provided
  the Novikov Conjecture is true for $\pi=\pi_1(M)$.
\end{cor}
\textit{Proof:} The uniformly finite cohomology of any leaf
$H_{uf}^k(L;\textbf{Q})=0$ for $k>m$ \cite{Att} and \cite{Gromov}, Example
0.1.C.
So following through the proof of Theorem 4.1
we have that $Pont^{k}(\pi)(TM/E)=0$ if $k>2m$, by Shulman's
theorem with $q=m$. Since the homotopy type of $M$ is that of $BG$ by
Haefliger's Corollary 3.2.4, the result follows.
\begin{cor} Let $M$ be a closed, oriented, smooth manifold with an integrable
  subbundle $E$ of the tangent bundle $TM$ of $M$, finitely generated
  fundamental group
  and a homomorphism $f:\pi_1(M)\to \textbf{Z}$ so that the corresponding
  rational cohmology class $f^*(j)$, $j$ the generator of $H^1(\textbf{Z})$,
  pulls back to 0 on each leaf of the foliation. More generally, suppose
  that $\gamma \in H^*(M;\textbf{Q})$ belongs to the cohomology algebra
  generated by $H^1(M;\textbf{Q})$ and pulls back to 0 on each leaf of the
  foliation.Then the higher signature
  class $f^*(j) \cup L_q(TM/E)$ or $\gamma \cup L_q(TM/E)$ vanishes for
  $q+1$ resp. $q+dim(\gamma)$ larger than $2dim(TM/E)$.
\end{cor}
\textit{Proof:}
Let $i:L \to M$ be a leaf of the foliation. We have a differential form
$\omega$ representing $f^*(j)$, which can be integrated to yield a map
$f:M \to S^1$, sending the leaves of the foliation to points. We have
$$L \overset{i}\to M \overset{f}\to S^1 \overset{j}\to S^1$$
The composition $j\circ f \circ i$ is null-homotopic. Since $j$ is
the generator, it is the identity map. So $f \circ i$ is null-homotopic.
We now foliate $S^1$ by points, and we have a map of foliations
$(M,\mathcal{F}) \to (S^1, \mathcal{C})$, i.e. a map which takes leaves
to leaves and is continuous along the leaves, where $\mathcal{C}$ is the
foliation by points of $S^1$. Passing to the leaf-spaces we get
a map $g:BG \to S^1$ induced from $f$. Let $k:M \to BG$ and
$k^\prime:S^1 \to S^1$ be the maps to the leaf spaces.
We consider the assembly map
$$A:H^*(BG; \mathcal{L}^{uf}(S_x)) \to H^*(M;\textbf{L})$$
where $\mathcal{L}^{uf}(S_\sigma)$ is the cosheaf over $BG$ which assigns the
uniformly finite \textbf{L}-cohomology of the leaves through the
simplex $\sigma$ to the
simplex $\sigma$, where $\textbf{L}$ is the simply-connected surgery spectrum.
and the analogous map for $S^1$:
$$A^\prime:H^*(S^1; \textbf{L}) \to H^*(S^1;\textbf{L})$$
We have a commutative diagram:
$$\begin{CD}
  H^*(BG; \mathcal{L}^{uf}(S_x)) @<g^*<< H^*(S^1,\textbf{L})\\
  @VAVV                   @VA^\prime VV\\
  H^*(M;\textbf{L})  @<f^*<< H^*(S^1;\textbf{L})
  \end{CD}$$
For every leaf $L$, $f \circ i$ is null-homotopic. So passing to KO-theory:
$$\Delta_{f^*(j)}(TM/E)= A(\sum_x\Delta^{(\infty)}_{g^*(j)}(\nu(S_x))x)$$
if we take the Pontrjagin character:
$$f^*(j)\cup L_i(TM/E)=A_*(\sum_x g^*(j)\cup L^{(\infty)}_i(\nu(S_x))x)$$
where $i=0,...,q$, $A_*$ is $A$ composed with the Pontrjagin character and
$H^*(BG;\mathcal{L})$ being where the right hand side now lies,
where $\mathcal{L}$ is the cosheaf assigning $H_{uff}^*(L_x;\textbf{Q})$ to
the leaf $L_x$ through $x \in BG$. Now the hypothesis on the leaves of the
foliation implies that $g^*(j)\cup L_i^{(\infty)}$ is zero on the coefficients
above dimension 0. So $f^*(j)=g^*(j) \in H^*(BG;\textbf{Q})$.

The Novikov Conjecture holds for $\textbf{Z}$ so that
the first result follows. To see the second result, note that
$H^1(M;\textbf{Z})=Hom(\pi_1(M),\textbf{Z})$, then note that Haefliger's
result applies to products and use the first result.

\begin{cor} Let $M$ be a closed  $K(\pi,1)$ (Eilenberg-MacLane) manifold.
  Suppose that $\pi$ satisfies
  the Novikov conjecture. If there is a foliation $\mathcal{F}$ of codimension
  $k$ with $2k<n$ on $M$, $n$ the dimension of $M$,
  so that the
  fundamental group of the leaf space $\pi_1(BG)$ satisfies the
  Novikov conjecture, then
  $\pi_1(L)\to \pi_1(M)$ is nontrivial for one of the the holonomy coverings of
  the leaves $L$ of $\mathcal{F}$, and  $H^n(B\pi_1(BG);\textbf{Q})=0$.
\end{cor}
\begin{rk} (See \cite{Gabai}). Haefliger \cite{Haef} generalized foliations to what are now
  known as Haefliger structures. Haefliger constructed a classifying space
  $B\Gamma_k^r$ for codimension $k$ $C^r$ foliations and showed that
  homotopy classes of such structures on a space $X$ are in 1-1 correspondance
  with $[X:B\Gamma_k^r]$, homotopy classes of maps of $X$ into $B\Gamma_k^r$.
  Codimension $k$ foliations $\mathcal{F}_0$ and $\mathcal{F}_1$ are
  concordant if there exists a codimension $k$ foliation on $M \times I$
  which restricts to $\mathcal{F}_0$ on $M \times 0$ and $\mathcal{F}_1$ on
  $M \times 1$. 
  Using the Haefliger's theory, Thurston proved \cite{Thurston1}
  that concordance classes of foliations on a manifold $M$ are in 1-1
  correspondance with homotopy classes of Haefliger structures $\mathcal{H}$
  together with concordance classes of maps $i:v_{\mathcal{H}}\to T(M)$.
  In addition he proved \cite{Thurston1} that if $K$ is a smooth $p$-plane
  field of codimension $q$, then $K$ is homotopic completely integrable
  $p$-plane field if and only if the map of $M$ to $BGL_q$ classifying the
  normal bundle to $K$ can be lifted to a map into $B\Gamma_q^r$ where
  $D:B\Gamma_q^r \to BGl_q$ is induced by the differential.
   \end{rk}
The present result for strengthening Bott's result on Pontrjagin classes
applies as well to Chern classes.
Let $M$ be a complex manifold $E$ a holomorphic vector bundle
which is integrable, and let $BG$ be the classifying space of the 
corresponding foliation groupoid with a homomorphism $j:\pi_1(BG) \to \pi$, 
$\pi$ a group.
We may consider the corresponding map $i:M \to B\pi$ and take
$i^*y$ for $y \in H^*(B\pi)$, and consider all products with Chern classes
$i^*(y)\cup c_i(E)$ in $H^*(M;\textbf{C})$, for all $y \in H^*(B\pi)$
 which we will call the higher Chern ring, and denote $Chern^*(\pi)(E)$
\begin{thm} Let $M^n$ be a complex manifold, $E$ an integrable holomorphic 
subbundle of $TM$, $\pi$ a group with a map $\pi_1(BG) \to \pi$, where $BG$ 
is the classifying space of the holonomy groupoid. If $\pi$ satisfies
the Novikov conjecture, hen the higher Chern ring
$Chern^q(\pi)(TM/E)$ vanishes for $q > 2k$. 
\end{thm}
The proof is the same as above, except with $\textbf{L}(\textbf{C})$, 
$K$-theory and the Chern character
are substituted for $KO$-theory and the Pontrjagin character.
\newtheorem{ex}{Example}[section]
\begin{ex} Let $M$ be a complex manifold of complex dimension $n$
  with a nonzero holomorphic vector field with fundamental group $\pi$,
  where $\pi$ satisfies the Novikov conjecture.
Then all the higher Chern numbers of $M$ vanish.
\end{ex}
\textit{Proof:} Let $E$ be the subbundle of $TM$ generated by the vector
field.
 Since $E$ is trivial, $Chern^{2n}(\pi)(TM/E)=
Chern^{2n}(\pi)(TM)=0$.
\begin{ex} Let $X$ be a manifold with the homotopy type of 
$S^4 \times S^4 \times S^4 \times S^4$ 
with Pontrjagin classes $p_1$,$p_2$ and $p_3$ nonzero. 
Consider the manifold $M=X \# T^{16}$. Let $E$ be a subbundle 
of $TM$ of codimension $4+i$. Consider the product of $x \cup p_1(TM/E)$ and 
$y \cup p_1^2(TM/E)$ in $Pont^*(\pi)(TM/E)$, where $x$ and $y$ are 
cohomology classes. This product is nonzero, when $x,y \in H^i(T^{16})$
have nonzero product hence the bundle is not integrable, that is has no
foliation with simply connected leaves.  
\end{ex}
\begin{ex} Consider a simply connected manifold $M$ for which Bott's 
obstruction prohibits integrability of a subbundle $E$ of $TM$
 (e.g., in the complex case all subbundles of $TM$ for $CP^n$ \cite{Bott})
and take the product $M \times T^k$ with
a $k$-torus. Then higher vanishing means that the corresponding subbundle of 
the tangent bundle of $M \times T^k$  is not integrable, i.e. has no foliation
whose leaves are simply connected. 
We similarly have this result for $M \times K(\pi, 1)$ where $K(\pi,1)$ is an 
aspherical manifold of fundamental group $\pi$. 
\end{ex}
\begin{ex} Consider the higher Chern classes for $M \times T^{2k}$,
where $E$ is a subbundle of $TM$ for which the Bott obstruction in
terms of Chern classes is nonzero, where $M$ is a complex manifold.
These will obstruct integrability. 
The corresponding subbundle of the tangent bundle of $M \times T^{2k}$
is not integrable, that is has no foliation whose leaves are simply connected.
\end{ex}
\begin{prop} Let $M^n$ be a closed $K(\pi,1)$ (Eilenberg-MacLane) manifold. Then if a subbundle
  $E \subset TM$ is integrable, the ring composed of the Pontrjagin
classes of $TM/E$ cup product with an element of the cohomology of $M$ is 
zero above dimension $2k$ where $k$ is the dimension of $TM/E$
and the leaves of the foliation are simply connected.
\end{prop}
\textit{Proof:} The cohomology of $M$ is the same as the cohomology of the
fundamental group. The result follows from Higher Integrability. 
\begin{prop} Let $M^n$ be a manifold with a map $M \to K(\pi,1)$ which 
induces an isomorphism on $\pi_i(M)$ for $i<m$. Then if a subbundle 
$E \subset TM$ is integrable, the ring composed of the Pontrjagin
classes of $TM/E$ cup product with the cohomology of $M$ up to dimension $m$
is zero above dimension $2k$ where $k$ is the dimension of $TM/E$. 
\end{prop}
\textit{Proof:} Same as above pulling the classes back from $K(\pi,1)$.
\begin{ex} Let $A$ be a $K(\pi,1)$ manifold of dimension $m$ and $B$
a $K(\pi,1)$ manifold of dimension $n$, $m<n$. Consider the foliation
of $M=A \times B$ with leaf $A$. Then higher integrability holds trivially. 
\end{ex}
\section{Foliated Surgery}
In this section we introduce surgery groups for foliated manifolds. The
groups defined in this section will lead to a general theory of foliated
surgery which is developed further in \cite{AttCap}.
\begin{dfn} Let $(M,\mathcal{F})$ be a compact manifold $M$ with foliation
  $\mathcal{F}$. Let $\mathcal{A}$ be an additive category, $\pi$ a group.
  Then the category $C_M^{\mathcal{F},bg}(\mathcal{A}[\pi])$ is defined
  to be the one whose objects are formal direct sums
  $$M=\sum_{x \in BG}\sum_{y \in L}M(y)$$
  of objects $M(y)$ in $\mathcal{A}[\pi]$, where $L$ is the leaf corresponding
  to the point $x$ in the classifying space of the holonomy groupoid $BG$
  which fall into a fixed finite number of types inside a ball of fixed
  radius in $L$ and in a ball of fixed radius in a transversal to $L$. Here
  $\mathcal{A}[\pi]$ is the category with the one object $M[\pi]$ for each
  object $M$ in $\mathcal{A}$, and with morphisms linear combinations of
  morphisms $f_g:M \to N$ in $\mathcal{A}$ finite. In the case $\pi=0$ we call
  $C_M^{\mathcal{F},bg}(\mathcal{A}[\pi])$ the category $C^{\mathcal{F},bg}$.
\end{dfn}
\begin{dfn} Let $(M,\mathcal{F})$ be a compact manifold $M$ with a foliation
  $\mathcal{F}$. Let $\mathcal{A}$ be an additive category, $\pi$ a group.
  Then the category $C_M^{\mathcal{F},bdd}(\mathcal{A}[\pi])$ is defined
  to be the one whose objects are formal sums
  $$M=\sum_{x \in BG}\sum_{y \in L}M(y)$$
  of object $M(y)$ in $\mathcal{A}[\pi]$, where $L$ is the leaf corresponding
  to the point $x$ in the classifying space of the holonomy groupoid $BG$.
  Here $\mathcal{A}[\pi]$ is the category with one object $M$ in $\mathcal{A}$,
  and with morphisms $f_g:M \to N$ in $\mathcal{A}$ finite. If $\pi=0$ we
  have $C_M^{\mathcal{F},bdd}(\mathcal{A}[\pi])=C_M^{\mathcal{F},bdd}$.
\end{dfn}
\begin{dfn}\cite{AM} Let $X_1$ and $X_2$ be spaces equipped with continuous maps
  $p_1,p_2$ to a metric space $Z$. Then a map $f:X_1 \to X_2$ is boundedly
  controlled if there exists an integer $m \ge 0$ so that for all
  $z \in Z$, $r\ge 0$, $p_1^{-1}(B_r(z)) \subseteq f(p_2^{-1}(B_{r+m}(z))$,
  where $B_r(z)$ denotes the metric ball in $Z$ of radius $r$ about $z$.
  Or equivalently, there is a constant $m \ge 0$ so that
  $$dist_Z(p_2 \circ f(x), p_1(x))< m$$
  for all $x \in X_1$.
\end{dfn}
\begin{dfn}(See \cite{AM}). Let $(X_1,\mathcal{F}_1)$ and
  $(X_2, \mathcal{F}_2)$ be foliated spaces, and let $(Z,\mathcal{G})$ be
  foliated space such that there are foliated maps $p_1:(X_1,\mathcal{F}_1)
  \to (Z, \mathcal{G})$ and $p_2:(X_2,\mathcal{F}_2)\to (Z,\mathcal{G}).$
  The foliated map $f:(X_1,\mathcal{F}_1) \to (X_2,\mathcal{F}_2)$ is
  boundedly controlled over the space $(Z,\mathcal{G})$ such that for
  any leaf $L$ and for any transversal $T$ in $X_1$, there is a constant
  $m \ge 0$ so that
  $$dist_Z(p_2 \circ f(x), p_1(x))< m$$
  for $x\in L$ or $x \in T$.
\end{dfn}
\begin{dfn} (See \cite{AM}). Let $(X,\mathcal{F})$ be a foliated space
  controlled over a foliated space $(Z, \mathcal{G})$ by a map $p$. Denote
  by $\mathcal{P}$ the category of metric balls in $Z$ with morphisms given
  by inclusions. Define $\mathcal{P}Hol(X,\mathcal{F})$ to be the category
  whose objects are pairs $(x,K)$ where $x \in X$ and
  $K \in \mid \mathcal{P}\mid$ is an
  object of $\mathcal{P}$ and if $x$ and $y$ are on the same leaf, a morphism
  $(x,K) \to (y,L)$ is a pair $(\omega, i)$ where $i \in \mathcal{P}(K,L)$
  is a morphism in $\mathcal{P}$ from $K$ to $L$ and $\omega$ is a
  holonomy on $p^{-1}(L)$ from $y$ to $p^{-1}(i(x)).$
\end{dfn}
\begin{dfn}(See \cite{AM}). The category of controlled free
  $\textbf{Z}\mathcal{P}Hol(X,\mathcal{F})$=modules is defined to be the
  category of controlled modules over the leaves of $\mathcal{F}$, with
  morphisms equal to the controlled morphisms. We denote this category by
  $\textbf{Z}\mathcal{P}Hol(X,\mathcal{F})^{bdd}.$
\end{dfn}
\begin{dfn} The category of controlled free $bg$
  $\textbf{Z}\mathcal{P}Hol(X,\mathcal{F})$-modules is defined to be the
  category of controlled modules so that in any ball of fixed radius of
  a leaf or transversal the modules fall into a finite number of types.
  Morphisms are defined to be morphisms of modules so that if the control
  space is partitioned into neighborhoods of fixed radius, the restrictions
  fall into a finite number of equivalence classes. We denote this
  category by $\textbf{Z}\mathcal{P}Hol(X,\mathcal{F})^{bg}$.
\end{dfn}
\begin{dfn} Let $(M,\mathcal{F})$ be a compact manifold with a foliation
  $\mathcal{F}$. The algebraic surgery group $L_*^{\mathcal{F},bg}$ is
  defined as
  $$L_*^{\mathcal{F},bg}(M)=L_*(\textbf{Z}\mathcal{P}Hol(M)^{bg}).$$
\end{dfn}
\begin{dfn} Let $(M,\mathcal{F})$ be a compact manifold with a foliation
  $\mathcal{F}$. The algebraic surgery group $L_*^{\mathcal{F},bdd}(M)$ is
  defined as
  $$L_*^{\mathcal{F},bdd}(M)=L_*(\textbf{Z}\mathcal{P}Hol(M)^{bdd}).$$
\end{dfn}
\begin{dfn} Let $(X,\mathcal{F})$ be a foliation on a compact manifold
  $X$ and $\mathcal{A}^{\mathcal{F},bg}(K)$ be the additive category of
  $K$-based objects in $\mathcal{A}^{\mathcal{F},bg}$ with morphisms
  $f:M \to N$ such that $f(\tau, \sigma)=0:M(\sigma) \to N(\tau)$ unless
  $\tau \ge \sigma$ so that $f(M(\sigma)) \subseteq \sum_{\tau \ge \sigma}
  N(\tau)$. Similarly for $\mathcal{A}^{\mathcal{F},bdd}$.
\end{dfn}
\begin{dfn} Let $(M,\mathcal{F})$ be a foliation on a compact manifold $M$.
  The quadratic foliated structure groups of $(R,K)$, where $R$ is a ring,
  are the cobordism groups for either $bg$ or $bdd$:
  $$\mathcal{S}_n^\mathcal{F}(R,K)=L_{n-1}(\mathcal{A}^\mathcal{F}(R,K),
  \mathcal{C}^\mathcal{F}(R,K),\mathcal{C}^{\mathcal{F}}(R)_*(K))$$
\end{dfn}
\begin{dfn} Let $(M, \mathcal{F})$ be a foliation on a compact manifold
  $M$. Define the local, finitely generated free $(R,K)$-modules for
  $bg$ or $bdd$:
  $$\Lambda(R)_*^\mathcal{F}(K)=(\mathcal{A}^\mathcal{F}(R,K), \mathcal{B}^\mathcal{F}(R,K), \mathcal{C}^\mathcal{F}(R)_*(K))$$
  where $\mathcal{B}^\mathcal{F}(R,K)$ is the category of finite chain
  complexes of f.g. free foliated $(R,K)$-module chain complexes $C$ such
  that $[C][\sigma](\sigma \in K)$ is a contractible f.g. free $R$-module
  chain complex.
\end{dfn}
\begin{dfn} Let $M$ be a compact manifold with foliation $\mathcal{F}$.
  Define the $bg$ foliated normal invariants for $(M, \mathcal{F})$ to be:
  $$H_n(BG; \mathcal{L}^{uf})=L_n(\Lambda(\textbf{Z})_*^{\mathcal{F},bg}(M))$$
  where $BG$ is the classifying space of the holonomy groupoid of $M$, and
  $\mathcal{L}^{uf}$ is the cosheaf assigning to each point $x$ of $BG$ the
  uniformly finite $L$-homology of the leaf $S_x$ through $x$.
\end{dfn}
\begin{dfn} Let $M$ be a compact manifold with foliation $\mathcal{F}$.
  Define the bdd foliated normal invariants for $(M,\mathcal{F})$ to be:
  $$H_n(BG;\mathcal{L}^{lf})=L_n(\Lambda(\textbf{Z})_*^{\mathcal{F},bdd}(M))$$
  where $BG$ is the classifying space of the holonomy groupoid of $M$ and
  $\mathcal{L}^{lf}$ is the cosheaf assigning to each point $x$ of $BG$ the
  locally finite $L$-homology of the leaf $S_x$ through $x$.
\end{dfn}
\begin{conj}[Bounded Novikov Conjecture for Foliations]
  Let $(M,\mathcal{F})$ be a foliation
  on a compact manifold. The map
  $$H_*(BG;\mathcal{L}^{lf})\to L_*^{\mathcal{F},bdd}(M)$$
  where $\mathcal{L}^{lf}$ is the cosheaf assigning the locally finite
  L-homology of the leaf through $x \in BG$ to the point $x \in BG$, is
  injective, provided the leaves of $(M,\mathcal{F})$ are uniformly
  contractible.
\end{conj}
\begin{conj}[Bounded Borel Conjecture for Foliations]
Let $(M,\mathcal{F})$ be a foliation on a compact manifold. The map
$$H_*(BG;\mathcal{L}^{lf}) \to L_*^{\mathcal{F},bdd}(M)$$
is an isomorphism, where $\mathcal{L}^{lf}$ is the cosheaf assigning the
locally finite L-homology of the leaf through $x \in BG$ to the point
$x \in BG$, provided the leaves of $(M,\mathcal{F})$ are aspherical.
\end{conj}
\begin{conj}[Bounded Geometry Novikov Conjecture for Foliations]
  Let $(M,\mathcal{F})$ be a foliation on a compact manifold. If the leaves
  of $(M,\mathcal{F})$ are aspherical, then
  $$H_*(BG;\mathcal{L}^{uf}(S_x)) \to L_*^{\mathcal{F},bg}(M)$$
  is injective, where $S_x$ is the leaf of $\mathcal{F}$ through $x \in BG$
  and $\mathcal{L}^{uf}$ is the uniformly finite L-homology cosheaf of the
  leaves over $BG$.
\end{conj}
\begin{conj}[Bounded Geometry Borel Conjecture for Foliations] Let
  $(M,\mathcal{F})$ be a foliation on a compact manifold. If the leaves of
  $(M,\mathcal{F})$ are aspherical, then
  $$H_*(BG;\mathcal{L}^{uf}(S_x))\to L_*^{\mathcal{F},bg}(M)$$
  is an isomorphism, where $S_x$ is the leaf of $\mathcal{F}$ through
  $x \in BG$ and $\mathcal{L}^{uf}$ is the uniformly finite L-homology
  cosheaf of the leaves over $BG$.
\end{conj}
\begin{dfn} Suppose $(M,\mathcal{F})$ is a foliation on a compact manifold
  with all leaves aspherical. Then $L_*^{\mathcal{F},bdd}(M)$ will be
  denoted $L_*(BG)$, where $BG$ is the classifying space of the holonomy
  groupoid of the foliation. The symmetric group $L^*(BG)$ is the
  symmetric L-group $L^{*,\mathcal{F},bdd}(M)$, and the quadratic
  L-group $L_*(BG)$ is $L_*^{\mathcal{F},bdd}(M)$.
\end{dfn}
\begin{dfn} Let $\textbf{L}(M,\mathcal{F})$ be the surgery space coming
  from the algebraic bordism category
  $$\Lambda(R)^{\mathcal{F},bdd}=(\mathcal{A}^{\mathcal{F},bdd}(\textbf{Z}),
  \mathcal{B}^{\mathcal{F},bdd}(\textbf{Z}), \mathcal{C}^{\mathcal{F},bdd}(
  \textbf{Z}))$$
  where $\mathcal{B}^{\mathcal{F},bdd}(\textbf{Z})$ is the category of finite
  chain complexes of f.g. free foliated $\textbf{Z}$-modules, and such that
  $\mathcal{C}^{\mathcal{F},bdd}$ is the category of contractible finite
  chain complexes of f.g. free foliated $\textbf{Z}$-modules. The space
  $\textbf{L}(M,\mathcal{F})$ will be called $\textbf{L}(BG)$ if the leaves
  of $(M,\mathcal{F})$ are aspherical, where $BG$ is the classifying space
  of the holonomy groupoid of $(M,\mathcal{F})$.
\end{dfn}
\section{Groupoid Integrability}
\begin{prop} The foliated Novikov conjecture holds for a classifying space
  for a holonomy groupoid $BG$ if and only if the assembly map in symmetric
  L-theory
  $$A_{BG}:H_*(BG;\textbf{L})\to L^*(BG)$$
  is a rational split injection. Note that $L^*(BG)\otimes \textbf{Z}[1/2]=
  L_*(BG)\otimes \textbf{Z}[1/2]$ where $L^*$ is the symmetric L-theory
  and $L_*$ is ordinary L-theory.
\end{prop}
\begin{dfn} Suppose $M$ is a manifold and $\mathcal{F}$ a foliation on $M$
  with the classifying space for the holonomy groupoid $BG$. Let
  $y\in H^i(BG;\textbf{Z})$. Then we have higher Pontrjagin classes
  $y \cup p_i(TM/E)$ where $E \subset TM$ is a subbundle of the tangent
  bundle of $M$. We have the higher Pontrjagin ring $Pont^*(BG)(TM/E)$
  consisting of all products with higher Pontrjagin classes $y \cup p_i(E)$
  for $y \in H^*(BG;\textbf{Z})$.
\end{dfn}
\begin{lmm}[Higher Shulman's Theorem] Let $(M,\mathcal{F})$ be a foliated
  manifold and let $BG$ be the leafspace of $\mathcal{F}$. If the higher
  Pontrjagin classes of $(M,\mathcal{F})$ vanish above the codimension $k$
  of $\mathcal{F}$, then the higher Pontrjagin ring
  $Pont(BG)^q(\nu(\mathcal{F}))$ vanishes above dimension $2k$, where
  $\nu(\mathcal{F})=TM/E$, $E$ being the tangent bundle of $\mathcal{F}$.
\end{lmm}
\par\noindent\textit{Proof:} We consider the cohomology of $BG$ with
coefficients in the cosheaf of uniformly finite $\textbf{L}$-cohomology
groups of the transversals, where $\textbf{L}$ is the surgery spectrum:
$$H^*(BG;\mathcal{L}_{uf}^*(T_x))$$
where $T_x$ is the transversal through $S_x$ at $x \in BG$ the leaf
through $x \in BG$ and $\mathcal{L}_{uf}^*$ is the $\textbf{L}$-cohomology
cosheaf. From the standard Leray-Serre spectral sequence we get an
assembly map
$$H^*(BG;\mathcal{L}_{uf}^*(T_x)) \to H^*(M;\textbf{L})$$
We then tensor with the reals:
$$H^*(BG;\mathcal{L}_{uf}^*(T_x))\otimes\textbf{R}\to H^*(M;\textbf{L})\otimes
\textbf{R}$$
We then use the fact that $\textbf{L}\otimes\textbf{Q}=\textbf{KO}\otimes
\textbf{Q}$ and we have
$$H^*(BG;\mathcal{KO}_{uf}^*(T_x))\otimes\textbf{R}\to H^*(M;\textbf{KO})
\otimes\textbf{R}$$
where $\mathcal{KO}_{uf}^*(T_x)$ is the cosheaf of uniformly finite
KO-cohomology groups of the transversals. We then apply the Pontrjagin
character to get
$$H^*(BG;\mathcal{H}_\beta^*(T_x))\otimes\textbf{R}\to H^*(M;\textbf{R})$$
where $\mathcal{H}_\beta^*(T_x)$ is the cosheaf of the bounded de Rham
cohomology of the transversals. We observe that the de Rham forms representing
elements of the $L^\infty$ Pontjagin ring of $T_x$ have a decomposition via
Shulman's theorem:
$$\Phi=\Phi_{q-1}+\Phi_{q-2}+...+\Phi_0$$
where $\Phi_i \in \Omega^{q+1}N_{q-i}SO(n)$ and $\Omega$ is the de Rham
functor and $\Omega NSO(n)$ is the de Rham complex of $NSO(n)$. Elements
of the higher Pontrjagin ring of $T_x$ can be represented by elements to
be represented as a de Rham form $\gamma \otimes \Phi$, where $\gamma$ is
a de Rham form representing an element of $H^*(BG)$ and $\Phi$ is a
representative of the Pontrjagin ring of $T_x$. But after assembling to
$M$ we see that this yields a de Rham form representing an element of the
higher Pontrjagin ring of $TM/E$ and that this is zero above dimension $2k$.
\begin{thm} Let $M$ be a compact closed manifold of dimension $n$, suppose
  $E \subset TM$ is a subbundle of the tangent bundle of $M$. Let $E$ be
  integrable, $BG$ the classifying space of the foliation groupoid of the
  corresponding foliation $\mathcal{F}$. Suppose the foliated Novikov
  conjecture is true for $\mathcal{F}$. Then the higher Pontrjagin ring
  $Pont^q(BG)(TM/E)$ vanishes above dimension $2k$, where $k$ is the
  dimnsion of $TM/E$.
\end{thm}
\begin{rk} Let $y_1,...,y_n \in H^*(BG)$ be classes for which the foliated
  Novikov Conjecture is known. Then if we let $K$ be the subset of $H^*(BG)$
  consisting of these classes, we have that $Pont^q(K)(TM/E)$ vanishes
  above dimension $2k$, where $k$ is the dimension of $TM/E$, and where
  $Pont^q(K)(TM/E)$ is the higher Pontrjagin ring consisting of all
  products of higher Pontrjagin classes $y\cup p_i(TM/E)$ for $y \in K$.
\end{rk}

\textit{Proof:} Let $k=dim(TM/E)$, $n=dim(TM)$.
Let $BG$ be the classifying space of the foliation groupoid $G$ of
$\mathcal{F}$ and consider the homology of $BG$ with coefficients
in the cosheaf $\mathcal{L}^{uf}(S_x)$ of the $\textbf{L}(BG)$-homology
groups of the leaves of $\mathcal{F}$: $H_*(BG;\mathcal{L}^{uf}(S_x))$.
The cosheaf $\mathcal{L}_*^{uf}(S_x)$ assigns to each simplex $\sigma$ of
$BG$, the group $H_*^{uf}(S_\sigma; \textbf{L}(BG))$, which is the uniformly
finite $\textbf{L}(BG)$-homology of the union of the leaves of $\mathcal{F}$
corresponding to $\sigma$ in $BG$. Taking the sum of the simplices of $BG$
with coefficients in $\mathcal{L}^{uf}(S_x)$, we get elements in
$H_*(M;\textbf{L}(BG))$. We let $\mathcal{A}$ be the additive category of free
additive category of free modules over the cosheaf $\mathcal{L}^{uf}$ and
the simplicial complex $K$ is the classifying space $BG$.

In this case $\mathcal{A}_*[K]$ is the functor that assigns to a simplex
$\sigma \in BG$ the group $\mathcal{L}^{uf}_*(S_\sigma)$, the
$\textbf{L}(BG)$-homology of the union of the leaves corresponding to
$\sigma$. We the have the category $\mathcal{A}_*(K)$ as the category of
direct sums
$$\sum_{\sigma \in BG}\mathcal{L}^{uf}_*(S_\sigma)$$
We then have the covariant assembly functor
$$\mathcal{B}(\mathcal{A}_*[K])\to \mathcal{B}((\mathcal{A})^*(K))$$
taken to $\Delta_*[K](\sigma)=S^{\mid\sigma\mid}\Delta(\sigma)$ and
therefore ends up in homology. The assembly map takes a chain complex
$C_*(BG;\mathcal{L}_{uf})=C_*^{uf}(BG;\mathcal{L}_{uf})$ (since $BG$ is
compact) to the complex
$$C_*[BG]_r=\sum_{\sigma \in BG}C_{r-\mid\sigma\mid}(\sigma; \mathcal{L}^{uf}(S_\sigma))$$
where the dimension shift $r-\mid\sigma\mid$ is the chain duality $T$. We
then use the Leray-Serre spectral sequence of a cosheaf,
$$H_*(BG; \mathcal{L}^{uf}(S_x))\Rightarrow H_*(M;\textbf{L}(BG))$$
which converges to $H_*^{uf}(M;\textbf{L}(BG))=H_*(M;\textbf{L}(BG))$ as
$M$ is compact and clearly yields the assembly map to $H_*(M;\textbf{L}(BG))$.

We get an assembly map
$$H_*(BG;\mathcal{L}^{uf}(S_x))\to H_*(M;\textbf{L}(BG))$$
where $\textbf{L}(BG)$ is the surgery spectrum of $BG$. We next define
the higher signatures of the leaves. Let $i:M\to BG$ and $j:S_x \to M$ be
the embedding of the leaf $S_x$ through $x \in BG$ so that for
$y \in H^*(BG)$, $\gamma=j^*i^*(y)\in H_{uf}^*(S_x)$. Let
$\Delta^{(\infty)}(\nu(S_x)) \in KO_*^{uf}(S_x)\otimes\textbf{Q}$ be the
$KO[1/2]$ orientation of the normal bundle $\nu(S_x)$ of the leaf
$S_x$ through $x, x\in BG$. This higher signature class is detected by
the surgery space $\textbf{L}(BG)$. This class lies in
$H_*^{uf}(S_x; H_*(BG;\textbf{L}))$ and by the foliated Novikov conjecture
this group injects into $H_*^{uf}(S_x;\textbf{L}_*(BG))$. This is because
the symmetric assembly map is a split injection on ordinary homology, and
hence we can use the Zeeman Comparison Theorem, which states that if the
$E_2$ terms on the summand to obtain an injection.

We have the following commutative diagram:
$$\begin{CD}
  H_*(BG;\mathcal{H}^{uf}(BG)_*(S_x)\otimes\textbf{Q})@<ph^{(\infty)}<< H_*(BG;\mathcal{KO}^{uf}(BG)(S_x))\otimes\textbf{Q}@>A^{(\infty)}_{Surgery}>> H_*(BG;\mathcal{L}^{uf}(S_x))\otimes\textbf{Q}\\
  @VA_{Leaves}^{\mathcal{H}}VV @VA_{Leaves}^{\mathcal{KO}}VV @VA_{Leaves}^{\mathcal{L}}VV\\
  H_*(M;H_*(BG;\textbf{Q})) @<ph<<H_*(M;\textbf{KO}_*(BG))\otimes \textbf{Q} @>A_{Surgery}>> H_*(M;\textbf{L}(BG))\otimes \textbf{Q}
\end{CD}$$
Here $ph$ is an isomorphism, and $A_{Surgery}:H_*(M;\textbf{KO}(BG))\otimes
\textbf{Q} \to H_*(M;\textbf{L}(BG))\otimes\textbf{Q}$ is injective as
$KO_*(BG)\otimes\textbf{Q}\to L(BG)\otimes\textbf{Q}$ is by the foliated
Novikov conjecture for $BG$.

In addition we have the cosheaf $\mathcal{L}^{uf}(S_x)$ as above, and the
other cosheaves $\mathcal{H}^{uf}$ which associates to $x$ the uniformly
finite homology of the leaf through $x$ and $\mathcal{KO}^{uf}(BG)_*(S_x)$
is the cosheaf which associates to $x$ the uniformly finite $KO(BG)$-homology
of the leaf through $x$. The map
$$A_{Surgery}^{(\infty)}:H_*(BG;\mathcal{KO}^{uf}(BG)(S_x))\otimes\textbf{Q}
\to H_*(BG;\mathcal{L}^{uf}(S_x))\otimes\textbf{Q}$$
is the assembly map coming from the injective map $KO(BG)\otimes\textbf{Q}
\to L(BG)\otimes\textbf{Q}$ and is itself injective. The maps $A_{Leaves}^{\mathcal{H}},A_{Leaves}^{\mathcal{KO}},A_{Leaves}^{\mathcal{L}}$ are assembly maps
of the cosheaves.

We have $\sum_x\Delta_\gamma^{(\infty)}(\nu(S_x))[x]$ as a cycle in generalized
homology. These each assemble to a higher signature class
$\Delta_\gamma(\nu(\mathcal{F}))$, with $\gamma \in H^*(BG)$. We then have
classes in $H_*(M;L(BG))$ coming from the leaves. We take the
$L^\infty$ Pontrjagin character $\Delta_\gamma^{(\infty)}(\nu(S_x))$,
and obtain the higher total L-class $\gamma \cup \mathcal{L}_{Tot}^{(\infty)}=
\gamma \cup (1+L_1^{(\infty)}(\nu(S_x))+L_2^{(\infty)}(\nu(S_x))+...)$, which
yields $\gamma \cup L_0^{(\infty)}(\nu(S_x))=\gamma,
\gamma\cup L_1^{(\infty)}(\nu(S_x)),...,\gamma \cup L_k^{(\infty)}(\nu(S_x))$.
Note that the Pontrjagin character is a rational isomorphism. For dimensional
reasons, these higher $L$-classes vanish above the dimension of $\nu(S_x)$
which is $k$.

Because of the Novikov conjecture, if the class of the surgery obstruction
$\Delta_\gamma^{(\infty),1}(\nu(S_x))\in H_*^{uf}(S_x;\textbf{L}(BG))$
coming from the higher signature 
\break $ph^{(\infty)}(\Delta_\gamma^{(\infty)}(\nu(S_x))\in H_*^{uf}(S_x;KO(BG))$
is zero then the higher signature is zero.
So the vanishing of the higher Pontrjagin classes from $KO(BG)$ for
dimensional reasons guarantees by the Novikov conjecture for foliations
that the higher Pontrjagin classes from the surgery group vanish in the
same range. These classes assemble via $A_{Leaves}^{\mathcal{L}}$ to
$\Delta_\gamma^1(\nu(\mathcal{F}))$ in $H_*(M;\textbf{L}(BG))$. Again we
have a higher signature class $\Delta_\gamma(\nu(\mathcal{F}))$ in
$H_*(M;KO(BG))$, which injects into $H_*(M;\textbf{L}(BG))$ via $A_{Surgery}$
because of the Novikov conjecture for $BG$.

The Pontrjagin character of $\Delta_\gamma(\nu(\mathcal{F}))$ is of the
form $ph(\beta \cap \Delta(\nu(\mathcal{F})))$ where $\Delta(\nu(\mathcal{F}))$
is the $KO[1/2]$-orientation of $\nu(\mathcal{F})$, and $\beta \in KO(M)$
is such $ph(\beta)=\gamma$, representing $\Delta_\gamma(\nu(\mathcal{F}))$ in
$H_*(M;KO(BG))$ and hence comes from an invariant polynomial in the curvature
form of $\nu(\mathcal{F})$.

Because of the vanishing over dimension $k$ of the $L^\infty$ higher Pontrjagin
classes of the leaves we have vanishing of the higher Pontrjagin classes of
$\nu(\mathcal{F})$ above dimension $k$, from the assembly map. We have
vanishing over dimension $2k$ of the higher Pontrjagin ring after applying
the Higher Shulman's Theorem (Lemma 6.1).

\begin{cor} Let $M$ be a smooth manifold and let $E$ be an integrable
  subbundle so that the corresponding foliation $\mathcal{F}$
  is a foliation where all the leaves
  of $\mathcal{F}$ are contractible. Then if $Pont^*(M)(TM/E)$ is the
  higher Pontrjagin ring generated by
  the higher Pontrjagin classes $y \cup p_i(TM/E)$, $y \in H^j(M)$ then
  $$Pont^q(M)(TM/E)=0$$
  for $q>2k$, where $k=dim(TM/E)$, provided the classifying space of the
  holonomy groupoid $BG$ satisfies the foliated Novikov conjecture.
\end{cor}

\begin{conj} Let $M$ be a smooth manifold and let $E$ be an integrable
  subbundle so that the corresponding foliation $\mathcal{F}$
  is a foliation where all the leaves of $\mathcal{F}$ are rationally
  contractible, i.e. if $L$ is a leaf of $\mathcal{F}$ then
  $\pi_i(L)\otimes\textbf{Q}=0$ for all $i$. Then if $Pont^*(M)(TM/E)$
  is the higher Pontrjagin ring generated by the higher Pontrjagin
  classes $y \cup p_i(TM/E)$, $y \in H^j(M)$ then
  $$Pont^q(M)(TM/E)=0$$
    for $q>2k$, where $k=dim(TM/E)$, provided the classifying space of the
  holonomy groupoid $BG$ satisfies the Novikov conjecture.
\end{conj}
The present result for strengthening Bott's result on Pontrjagin classes
applies as well to Chern classes. Let $M$ be a complex manifold, $E$ a
holomorphic vector bundle which is integrable, and let $BG$ be the classifying
space of the corresponding foliation groupoid. We may consider the
corresponding map $i: M\to BG$ and take $i^*y$ for $y \in H^*(BG)$, and
consider all products with Chern classes $i^*(y)\cap c_i(E)$ in
$H^*(M;\textbf{C})$ for all $y \in H^*(BG)$ which we will call the higher
Chern ring, and denote $Chern^*(BG)(E)$.
\begin{thm} Let $M^n$ be a complex manifold, $E$ an integrable holomorphic
  subbundle of $TM$, $BG$ the classifying space of the holonomy groupoid
  of the corresponding foliation. If $BG$ satisfies the foliated Novikov
  conjecture, then the higher Chern ring $Chern^q(BG)(TM/E)$ vanishes for
  $q>2k$
\end{thm}
The proof is the same as above, except with $\textbf{L}(\textbf{C})$,
$K$-theory and the Chern character are substituted for $KO$-theory and the
Pontrjagin character.
\begin{thm} Let $M$ be a manifold of complex dimension $n$ with a nonzero
  holomorphic vector field with a holomorphic foliation with holonomy
  groupoid $BG$, where $BG$ satisfies the foliated Novikov Conjecture.
  Then all the higher Chern numbers of $M$ vanish.
\end{thm}
\textit{Proof:} Let $E$ be the subbundle of $TM$ generated by the vector
field. Since $E$ is trivial, $Chern^{2n}(BG)(TM/E)=Chern^{2n}(BG)(TM)=0$.
\begin{rk} We do not need $M$ to have infinite fundamental group to have
  nonzero integrability obstruction. Consider the Reeb foliation on
  $S^3$. This foliation has leaf space $BG=S^3$, and hence there
  is a class $x \in H^3(BG)$ which can be used to yield an integrability
  obstruction for $E \subset TM$.
\end{rk}
\begin{ex} Let $X$ be a manifold which is homotopy equivalent to
  $S^4 \times S^4 \times S^4 \times S^4$ with Pontrjagin classes
  $p_1$, $p_2$ and $p_3$ nonzero. Suppose $E \subset TM$ be a subbundle
  of the tangent bundle of $M=X$ of codimension 4. Then if $BG=S^4 \times
  S^4 \times S^4$ is the leafspace, and we choose $x,y \in H^4(BG)$, with
  nonzero product, so
  we have that the product of $x \cup p_1(TM/E)$ and $y \cup p_1(TM/E)$
  is nonzero and of dimension greater than 8, so $E$ is not integrable
  with leafspace $BG$.
  \end{ex}
\section{Higher Godbillon-Vey Invariant}
This section utilizes \cite{CS} in describing the Cheeger-Simons and 
Godbillon-Vey invariants. For a review of the Godbillon-Vey invariant,
see \cite{Ghys,Roger,Hurder1}. 
Let $G$ be a Lie group with finitely many components, $BG$ its classifying
space and $I^*(G)$ the ring of invariant polynomials on $G$. 
The Weil homomorphism constructs a homomorphism $w:I^k(G) \to H^{2k}(BG;
\textbf{R})$. Let $G=GL(n,\textbf{R})$ and set $I_0(G)=ker(w)$. Then 
$I_0=\sum I_0^k$ is the ideal generated by the polynomials $trA^{2k-1}$.
Let $\Lambda \subset \textbf{R}$ be a proper subring of the reals. Let
$M$ be a $C^\infty$ manifold and let $\Lambda^*$ denote the ring of 
differential forms on $M$. Let $C_k \supset Z_k \supset B_k$ denote the rings
of smooth singular chains, cycles and boundaries, and $\partial:C_k \to 
C_{k-1}$ and $\delta:C^k \to C^{k+1}$ be usual boundary and coboundary
operators. Let $\Lambda^k_0$ be the closed $k$-forms with periods lying in
$\Lambda$.
\begin{dfn}
$$\hat{H}^k(M; \textbf{R}/\Lambda)=\{f \in Hom(Z_k, \textbf{R}/\Lambda) \mid
f \circ \partial \in \Lambda^{k+1}\}$$
\end{dfn}
\begin{dfn}
$$R^k(M;\Lambda)=\{(\omega, u)\in \Lambda^k_0 \times H^k(M;\Lambda) \mid
r(u)=[\omega]\}$$
where $r$ is the natural map $r:H^k(M;\Lambda) \to H^k(M; \textbf{R})$ and 
$[\omega]$ is the de Rham class of $\omega$.
\end{dfn}
\begin{dfn} Let $f \in \hat{H}^k(M;\textbf{R}/\Lambda)$. Since 
$\textbf{R}$ is divisible, there is a real
cochain $T$ with $T \mid Z_k=f$. Since $\tilde{\delta T}=\delta\tilde{T}=
f \circ \partial$, by assumption there exists $\omega \in \Lambda^{k+1}$
and $c \in C^{k+1}(M; \Lambda)$ such that $\delta T=\omega-c$. Then 
$0=\delta^2 T=\delta \omega-\delta c=d\omega-\delta c$. Since a nonvanishing
differential form never takes values lying only in a proper subgroup
$\Lambda \subset \textbf{R}$, we conclude that $d\omega=\delta c=0$.
Since $\delta T=\omega-c$ it follows that $\omega \in \Lambda_0^{k+1}$, 
$[c]=u \in H^{k+1}(M;\Lambda)$ and $[\omega]=r(u)$. Set $\delta_1(f)=\omega$
and $\delta_2(f)=u$.
\end{dfn} 

Let $\alpha=\{E,M,\theta\}$ be a principal $G$-bundle total space $E$, base
space $M$ and connection $\theta$. Let $\epsilon(G)$ be the category of 
these objects with morphisms being connection preserving bundle maps.
$W:I^k(G) \to \Lambda^{2k}(M)$ is the natural transformation associated
with the Weil homomorphism. If $P \in I^k(G)$, $u \in H^*(BG; \Lambda)$ and 
$\Omega$ the curvature form of $\alpha \in \epsilon(G)$, then 
$W(P)=P(\Omega,...,\Omega)$ and $C_\Lambda(u)=u(\alpha)$ the characteristic 
class. Let $r:H^*(M;\Lambda) \to H^*(M;\textbf{R})$. Set
$$K^{2k}(G,\Lambda)=\{(P,u) \in I^k(G) \times H^{2k}(BG,\Lambda) \mid
w(P)=r(u)\}$$
\begin{thm}[Cheeger-Simons] \cite{CS} Let $(P,u) \in K^{2k}(G,\Lambda)$. 
For each $\alpha \in \epsilon(G)$ 
there exists a unique $S_{P,u} \in \hat{H}^{2k-1}(M; \textbf{R}/\Lambda)$
satisfying
\par\indent 1) $\delta_1(S_{P,u})=P(\Omega)$
\par\indent 2) $\delta_2(S_{P,u})=u(\alpha)$
\par\indent 3 If $\beta \in \epsilon(G)$ and $\phi:\alpha \to \beta$ is a 
morphism, then $\phi^*(S_{P,u}(\beta))=S_{P,u}(\alpha)$.
\end{thm}
Taking $\Lambda=\{0\}$ we see that $Q \to (Q,0)$ is an isomorphism 
between $I_0^k$ and $K^{2k}(G,\{0\})$. Let $\pi_1(M)\ne 0$ be the 
fundamental group of $M$ and let $i:\pi_1(M) \to \pi$ be a homomorphism
to a group $\pi$ inducing a map $j:M \to B\pi$. Let $\gamma \in H^*(B\pi)$
or $\gamma \in H^*(B\mathcal{G})$, where $B\mathcal{G}$ is the leafspace of the foliation
on $M$ be an element. If $\alpha \in \epsilon(GL(n,
\textbf{R}))$ and $Q \in I_0^k$ set let $j^*(\gamma) \cup Q$ be the 
function $\Omega \mapsto j^*(\gamma) \cup Q(\Omega)$, where $\Omega$
is a curvature form. Then we have
$$\hat{Q}_\gamma(\alpha)=S_{j^*(\gamma) \cup Q,0}(\alpha) \in \hat{H}^{2k-1}(M; \textbf{R})$$
Let $V=\{V, M, \nabla\}$ be a real vector bundle with connection. In
$BGL(n,\textbf{R})$ we have the Pontrjagin class $p$ and the polynomial
$P_k \in I^{2k}(GL(n,\textbf{R}))$ with $w(P_k)=p_k$. Letting 
$E(V) \in \epsilon(GL(n,\textbf(R)))$ be the basis bundle of $V$ we define
as before $j^*(\gamma)\cup P_k$ to be the function $\Omega \mapsto j^*(\gamma)
\cup P_k(\Omega)$ so that we have
$$\hat{p_k}_\gamma(V)=S_{j*(\gamma)\cup P_k, p_k}(E(V))$$
Let $\mathcal{F}$ be a foliation of codimension $n$ in $M$,  and let
$\nu(\mathcal{F})$ be its normal bundle. Let $N(\mathcal{F})=\{\nu(\mathcal{F}), M, \nabla\}$, where $\nabla$ is a Bott connection \cite{Bott} and set
$\hat{Q}_\gamma(\mathcal{F})=\hat{Q}_\gamma(N(\mathcal{F}))$, and 
$\hat{p_i}_\gamma(\mathcal{F})=\hat{p_i}_\gamma(N(\mathcal{F}))$. 
$$\hat{Q}_\gamma(\mathcal{F}) \in H^{2k-1}(M; \textbf{R}), Q \in I_0^k,
k>n$$
$$\hat{p_i}_\gamma(\mathcal{F}) \in H^{4i-1}(M, \textbf{R}/\textbf{Z}), 2i > n$$
\begin{thm}
These classes are defined independently of the choice of Bott connection
and are invariants of $\mathcal{F}$.
\end{thm}
\textit{Proof:} This follows directly from the higher integrability theorem.

For example for $Q=trA^2$, $\hat{Q}_\gamma(\mathcal{F})
 \in H^3(M;\textbf{R})$ is our higher Godbillon-Vey invariant.
\begin{ex}
We take the following discussion from \cite{BottH}.
Let $X$ be a manifold and $\pi:E \to X$ a vector bundle over $X$ with discrete
structural group $G$. Locally $E$ is a product and each product neighborhood
has a natural horizontal foliation. Since $G$ is discrete, the change of 
coordinates maps are locally constant and so preserve these local foliations. 
This gives a foliation of $E$ called the horizontal foliation. The quotient
bundle $\eta$ to the horizontal foliation of $E$ is isomoorphic $\pi^*(E)$.
Since $\pi^*$ is an isomorphism in cohomology, the Pontrjagin ring of 
$\pi^*(E)$ is the same as that of $E$. 

Let $X=S^{2n-1}/\textbf{Z}_p$ where $p$ is an odd prime $S^{2n-1}\subset
\textbf{C}^n$ and $\textbf{Z}_p$ acts through $U_1 \times ... \times U_1$.
Then $L$ be the basic complex line bundle over $X$. $L$ has structural 
group $\textbf{Z}_p$. The first Chern of $L$, $c_1(L)$ generates $H^2(X;
\textbf{Z}_p)$ and an easy argument using the Gysin sequence shows
$(c_1(L))^{n-1} \ne 0$ and 
$H^*(X;\textbf{Z}_p)\supset \textbf{Z}_p[c_1]/((c_1)^n)$.
The Pontrjagin class of $L$ is $1-(c_1)^2$ and the ring generated by it has
nonzero elements in dimension $2[n/2]$. Now take the manifold $M(\pi)$ or
$M(B\mathcal{G})$ of 
dimension $2n-1$ and fundamental group $\pi$ satisfying the Novikov conjecture
or $B\mathcal{G}$ with foliation $\mathcal{F}$ satisfying the foliated Novikov conjecture and multiply it by $X$ and consider
$p_i(N) \cup f^*(\gamma) \in H^*(X \times M(\pi))$, where $N$ is the normal 
bundle to the horizontal foliation and $f:M(\pi)\to B\pi$ is a map and 
$\gamma \in H^*(B\pi)$ or else $p_i(N)\cup f^*(\gamma) \in H^*(X \times M(B\mathcal{G}))$, where $N$ is the normal bundle to the horizontal foliation and
$f:M(B\mathcal{G})\to B\mathcal{G}$ is a map and $\gamma \in H^*(B\mathcal{G})$.
This is zero in real cohomology (but not in integral
cohomology), and so we let $B$ be the Bockstein, and 
$B({\hat{p_i}_\gamma})=-\psi^*(p_i \cup \gamma)\ne 0$ 
for this manifold, where $\psi:B\Gamma_n \to BGl(n,\textbf{R})$ is the natural 
map. 
\end{ex}
\begin{ex}
To find a nonzero higher Godbillon-Vey invariant, take a foliation on $S^3$
that has nonzero Godbillon-Vey invariant \cite{Th}, and take the product 
foliation on $S^3 \times M(\pi)$ or $S^3 \times M(B\mathcal{G})$,
where $M(\pi)$ and $M(B\mathcal{G})$ are the manifolds in Example
7.1. Then $\hat{Q}_\gamma \ne 0$ as above.
\end{ex}
\begin{ex}[Higher Maslov Index] Recall the consruction of the Maslov index from \cite{KT1}.
  We generalize this to a higher Maslov index.
  Let $G$ be a Lie group, $H \subset G$ a closed subgroup with
  finitely many connected components and $\Gamma \subset G$ a discrete
  subgroup, operating properly discontinuously and without fixed points on
  $G/H$, so that $M=\Gamma\backslash G/H$ is a manifold. The foliation of
  $M$ consists of one single leaf equal to $M$. The bundle
  $$P=\Gamma\backslash G \times_H G\simeq G/H\times_\Gamma G=M=\Gamma\backslash
  G/H$$
  is flat. Let $G=U(m)$ and $H=O(m)$. Take $N(\pi)$ to be a manifold with
  fundamental group $\pi$. Suppose $\pi$ satisfies the Novikov Conjecture.
  We have for $\gamma \in H^i(B\pi)$ 
  a higher Maslov class transgressing to $\gamma \cup c_1$ the higher odd Chern
  class, yielding $\mu_\gamma \in H^{i+1}(M \times N(\pi))$.
  This can also be generalized using the holonomy groupoid.
\end{ex}
\begin{rk} One can define a higher Chern-Simons class following the
  construction of Chern-Simons classes from \cite{Hurder2}. Let $M$ be
  a closed oriented 4-manifold and consider $M$ as foliated by points,
  so that this is a Riemannian foliation of codimension 4. Then if
  $\gamma \in H^1(B\pi_1(M))$ we can use an oriented framing of
  $M \times \textbf{R}^4$, to get a transgression of the higher first
  Pontrjagin class $\gamma \cup p_1$, the higher Chern-Simons class
  $cs_\gamma \in H^4(M)$.
  \end{rk}
We remark that the assembly construction of the present paper is the beginning
of a surgery theory for foliations, which is in analogy to Connes' index
theory \cite{Connes}, that will be discussed in a forthcoming paper
\cite{AttCap}.

\end{document}